\documentclass[leqno,12pt]{amsart}
\usepackage{amssymb,amsmath,amscd}
\usepackage{a4wide}

\newtheorem{thm}{Theorem}[section]

\newtheorem{corollary}[thm]{Corollary}
\newtheorem{prop}[thm]{Proposition}
\newtheorem{lemma}[thm]{Lemma}
\newtheorem{fact}[thm]{Fact}

\theoremstyle{definition}
\newtheorem{defn}[thm]{Definition}
\newtheorem{example}[thm]{Example}

\theoremstyle{remark}
\newtheorem{remark}[thm]{Remark}

\newcommand{\bt}{\begin{thm}}
\newcommand{\et}{\end{thm}}
\newcommand{\bp}{\begin{prop}}
\newcommand{\ep}{\end{prop}}
\newcommand{\bd}{\begin{defn}}
\newcommand{\ed}{\end{defn}}
\newcommand{\bl}{\begin{lemma}}
\newcommand{\el}{\end{lemma}}
\newcommand{\bfa}{\begin{fact}}
\newcommand{\efa}{\end{fact}}
\newcommand{\bc}{\begin{corollary}}
\newcommand{\ec}{\end{corollary}}
\newcommand{\bex}{\begin{example}}
\newcommand{\eex}{\end{example}}
\newcommand{\br}{\begin{remark}}
\newcommand{\er}{\end{remark}}
\newcommand{\ben}{\begin{enumerate}}
\newcommand{\een}{\end{enumerate}}
\newcommand{\ds}{\displaystyle}
\newcommand{\codim}{\mbox{codim}}

\newcommand{\Hom}{{\mathcal H}om}
\newcommand{\rmhom}{{\mbox{Hom}}}
\newcommand{\homr}[2]{\rmhom_{R} ( #1, #2)}
\newcommand{\zhomr}[2]{\ _0\rmhom_{R} ( #1, #2)}
\newcommand{\zhomx}[2]{\ _0\rmhom_{R_X} ( #1, #2)}

\newcommand{\ext}{{\mbox{Ext}}}

\newcommand{\sotto}[2]{#1_{#2}}
\newcommand{\lra}{\longrightarrow}
\newcommand{\rrr}{\rightarrow}
\newcommand{\ra}{\rightarrow}

\newcommand{\ii}{{\mathcal I}}
\newcommand{\ideal}[1]{\sotto {{\mathcal I}}{#1}}

\newcommand{\exact}[3]
{0 \rrr #1 \rrr #2
\rrr #3 \rrr 0}

\newcommand{\PP}{\mathbb{P}}
\newcommand{\pp}{\mathbb{P}}

\newcommand{\Z}{\mathbb{Z}}

\newcommand{\N}{\mathbb{N}}

\newcommand{\cnn}{{\mathcal N}}
\newcommand{\coo}{{\mathcal O}}

\newcommand{\caf}{{\mathcal F}}

\begin{document}

\title{Codimension $3$ Arithmetically Gorenstein Subschemes of projective $N$-space}
\author{R. Hartshorne, I. Sabadini and E. Schlesinger}
\date{}

\address{Department of Mathematics, University  of
California, Berkeley, California 94720--3840}

\address{Dipartimento di Matematica, Politecnico di Milano, Piazza Leonardo da
Vinci 32, 20133 Milano, Italia}

\address{Dipartimento di Matematica, Politecnico di Milano, Piazza Leonardo da
Vinci 32, 20133 Milano, Italia}

\thanks{
The second author was partially supported by
MIUR PRIN 2004 {\em Geometria  sulle variet\`a algebriche}.
The third author was partially supported by
MIUR PRIN 2005:
  {\em Spazi di moduli e teoria di Lie}.
 }


\begin{abstract}
We study the lowest dimensional open case of the question whether every arithmetically
Cohen--Macaulay
subscheme of $\PP^N$ is glicci, that is, whether every zero-scheme in $\PP^3$ is
glicci. We show that a general set of $n \geq 56$ points in $\PP^3$
admits no strictly descending Gorenstein liaison or biliaison. In order to prove
this theorem, we establish a number of important results about arithmetically Gorenstein
zero-schemes in $\PP^3$.

\end{abstract}

\maketitle

\section{Introduction}
There has been considerable interest recently in the notion of
Gorenstein liaison for subschemes of $\mathbb{P}^N$. In
particular, a question that has attracted a lot of attention, by
analogy with the case of complete intersection liaison in
codimension $2$, is whether every ACM subscheme of $\mathbb{P}^N$
is in the Gorenstein liaison class of a complete intersection
(glicci for short). Many special classes of ACM subschemes have
been found that are glicci, but the question in general remains
open \cite{CDH},\cite{CH},\cite{H3},\cite{KMMNP},\cite{M},\cite{MN}.

Our motivation for this research was to consider the lowest dimensional case of this question,
namely zero-dimensional subschemes of $\mathbb{P}^3$. In this case, since every zero-dimensional
scheme is ACM, the question becomes simply: is every zero-scheme in $\mathbb{P}^3$ glicci?
One of us has shown in an earlier paper \cite{H3} that a general set of $n\leq 19$ points
is glicci (we give a new proof of this in Proposition \ref{glicci}),
while for $n \geq 20$ it is unknown. In this paper we show

\bt[Theorems \ref{dl} and \ref{db}]
A general set of $n \geq 56$ points  in $\PP^3$ admits no strictly descending Gorenstein liaison
or biliaison.
\et

The theorem does not imply that a general set of $n \geq 56$ points is not glicci,
but it shows that,
if such a set is glicci, it must first be linked upwards to a larger set before eventually linking
down to a point.

In order to prove the theorem, since Gorenstein liaisons are performed with arithmetically Gorenstein (AG) subschemes,
we need to establish a number of results about these.

It is known that the family $\mbox{PGor}(h)$ of codimension three AG  subschemes of $\PP^N$ with a fixed Hilbert function,
encoded in the $h$-vector $h$ (see section 2), is irreducible \cite{D}. The dimension of $\mbox{PGor}(h)$ is computed in
\cite{K} and \cite{KM}, but not explicitly as a function of the $h$-vector.
In the case of zero-dimensional
AG subschemes of $\PP^3$, a formula for $\dim  \mbox{PGor}(h)$ in terms of the $h$-vector is given in
\cite[\S 5]{CV}; following a different approach, we derive in
Corollary \ref{dimension} a formula which allows
to compute $\dim  \mbox{PGor}(h)$ inductively.

In all cases we are aware of where a class of ACM subschemes of
$\mathbb{P}^N$ has been proved to be glicci, the proof was actually accomplished using strict
Gorenstein liaisons, i.e. using only  those AG schemes of the form $mH_X-K_X$ on some ACM scheme $X$,
where $H_X$ and $K_X$ denote respectively the hyperplane and the canonical divisor class.
So we ask whether all AG subschemes can be obtained in this way, or only some. This analysis was performed in an earlier
paper \cite{hartshorne1} for curves
in $\mathbb{P}^4$, and we extend the results to any $\mathbb{P}^N$, $N\geq 3$.
One may consult \cite{Kl} for a deformation theoretic approach to this problem, and Boij in \cite[Theorem 3.4]{B} proves a
result related to ours in an algebraic context.
To give a precise statement, we use
two  numerical invariants $b(h)$ and $\overline{b}(h)$ of the Hilbert function, which satisfy $b(h) \geq 2 \overline{b}(h)$.
Our result is
\bt [Theorem \ref{thm25}]
Given an $h$-vector arising as the $h$-vector of a codimension $3$ AG subscheme of $\PP^N$, there is always
in $\mbox{\em PGor}(h)$ a scheme of the form $mH_X-K_X$ on some codimension $2$ {\em ACM} subscheme $X \subset \PP^N$. Furthermore, if
$b(h) \geq 2 \overline{b}(h)+2$,
the general element of $\mbox{\em PGor}(h)$ has this form.
\et
It looks as if the most interesting $h$-vectors will be those with $b(h)$ equal to $2\overline{b}(h)$ or $2\overline{b}(h)+1$,
in which case it is possible the most general element of  $\mbox{PGor}(h)$ is not of the form $mH_X-K_X$ on any ACM $X$. For points in $\PP^3$,
we give strong evidence in Proposition \ref{30} that this happens in degree $30$, by giving an example of a family of AG sets of $30$ points
whose general element  is not of the form $mH_C-K_C$ on any integral ACM curve $C$.
This phenomenon may occur earlier in degrees $20$, $28$ or $29$.

Watanabe showed that codimension $3$ AG subschemes are licci, that is, in the complete intersection liaison class of a complete intersection.
In section 5 we sharpen this result in the case of general AG sets of points in  ${\mathbb P}^3$ by showing that
\bt[Theorem \ref{thm-bil}]
For a given $h$-vector, a general AG zero dimensional subscheme of ${\mathbb P}^3$ with that $h$-vector can be obtained by
ascending complete intersection biliaisons from a point.
\et
This was proved for curves in ${\mathbb P}^4$ in \cite{hartshorne1}. To establish the result,
we need a version of the well known Cayley--Bacharach property for a zero-scheme
in $\mathbb{P}^N$ to be AG, of which we give a new proof in section 4.

In section 6, we study the Hilbert schemes $\mbox{PGor}(h)$ of $AG$ zero-schemes $Z$ in $\PP^3$
keeping in mind the two crucial questions
\begin{enumerate}
\item When does the general element of $\mbox{PGor}(h)$ have the form $mH-K$ on
some ACM curve $C$?

\item  As the AG scheme $Z$ varies in the Hilbert
scheme $\mbox{PGor}(h)$, how many general points can we assign to a general $Z$ ?
This is important in order to understand the possible Gorenstein liaisons
one can perform on a set of general points.
\end{enumerate}
We don't have a general answer to these questions, but we have made computations up to degree $30$,
and summarized our results in Table 1. The notations and methods used for these computations are explained in section 6.

We hope that this ground work will be the foundation for an eventual solution of our motivating
problems discussed above.

We would like to thank the referee for his careful reading of the paper and his several helpful comments.

\section{The $h$-vector of ACM subschemes of ${\mathbb P}^N$}
\label{sec1}
Notation: $R= \mbox{H}^0_* (\PP^N,\coo_{\PP^N}) = \bigoplus_{n \in \Z} \mbox{H}^0 (\PP^N, \coo_{\PP^N}(n))$ denotes the homogeneous
coordinate ring of $\PP^N$; $R_X = R/I_X$ the coordinate ring of a closed subscheme
$X$, where $I_X = \mbox{H}^0_* (\PP^N, \ideal{X,\PP^N})$.

Throughout the paper, $X$ denotes an arithmetically Cohen-Macaulay (ACM for short)
subscheme of $\PP^N$: recall
$X$ is called ACM if the coordinate ring $R_X$ is Cohen-Macaulay (of dimension $t+1$ where $t= \dim X$).
We denote by $\Omega_X$ the graded canonical module of $R_X$. When $\dim X > 0$, we have
$$\Omega_X \cong \mbox{H}^0_* (\PP^N, \omega_X)= \bigoplus_{n \in \Z} \mbox{H}^0 (\PP^N, \omega_X (n))$$
where $\omega_X$ is the Grothendieck dualizing sheaf of $X$.

A subscheme $Z \subseteq {\mathbb P}^N$ is {\em arithmetically
Gorenstein} (AG for short) if its homogeneous coordinate ring is a Gorenstein ring.  This
is equivalent to saying $Z$ is ACM and the canonical module $\Omega_Z$
is isomorphic to $R_Z(m)$ for some $m \in {\mathbb Z}$.

Denote by $H(n)=\dim_k (R_X)_n$ the Hilbert function of $R_X$. The
difference function $h_X (n)=\partial^{t+1} H(n)$ is called the
$h$-{\em vector} of $X$ \cite[\S 1.4]{M}:  it is nonnegative and
with finite support. We let $b=b(X)$ denote the largest integer
$n$ such that $h_X (n) > 0$. One can show $h_X (n) > 0$ for
$0 \leq n \leq b$. It is convenient to represent the $h$-vector in the form
\begin{equation}\label{hvettore}
    h_X = \{ 1=h_X(0),h_X (1), \ldots, h_X (b) \}.
\end{equation}

We now recall how other numerical invariants can be computed in terms of the $h$-vector.
First of all, the degree of $X$ is given by the formula $d=\sum_{n=0}^b h_X(n)$.
If $X$ is nondegenerate, $h_X(1)={\rm codim} X$.

For any closed subscheme $X \subset \PP^N$, we denote by $s(X)$ the least degree of a hypersurface
containing $X$. If $X$ is ACM of codimension $c$, the number $s(X)$ is the least positive integer $n$ such that
$h_X (n) < \binom{n+c-1}{c-1}$.

The integer $b(X)$ is related to the regularity and to the index of speciality of $X$.
A sheaf $\caf$ on $\PP^N$ is $n$-regular in the sense of
Castelnuovo-Mumford if $$H^i (\PP^N, \caf(n-i)) = 0 \;\;\;\;\mbox{for $i>0$}.$$
If $\caf$
is $n$-regular, then it is also $n+1$-regular. Thus one defines the
regularity of a sheaf $\caf$ as the least integer $r$ such that $\caf$
is $r$-regular. The regularity of a subscheme $X \subseteq \PP^N$ is the regularity
of its ideal sheaf.

We let $m(X)$ denote the largest integer $n$
such that $(\Omega_X)_{-n} \neq 0$ (when $\dim X > 0$, this is the index of speciality
of $X$).

\bp \label{bb} Let $X$ be an ACM subscheme of $\PP^N$ of dimension
$t<N$. Then
\begin{enumerate}
  \item the regularity of $X$ is $b(X)+1$.
  \item we have $m(X)= b(X) - t - 1$.
  \item we have $m(X)=\max\{ n\ :\  h^{t+1} (\PP^N, \ii_X (n) )\not=0\}.$
\end{enumerate}
\ep
\begin{proof}
The first statement \cite[p.8 and p.30]{M} follows computing both the $h$-vector
and the regularity of $X$ out of the minimal free resolution of $I_X$ over $R$.
The second statement follows from the isomorphism
$$
\Omega_X \cong \mbox{Hom}_{R_L} (R_X, \Omega_L)
$$
where $L \cong \PP^{t}$, the map $R_L \ra R_X$ is induced by
projection from a general linear subspace of codimension $t+1$, so that
$R_X$ is a free finitely generated graded $R_L$ module with $h_X(n)$ minimal generators
in degree $n$.

The third statement is immediate when $X$ is a hypersurface. If
$0\leq t\leq N-2$ it is a consequence of the fact that the
canonical module $\Omega_X$ is the graded $k$-dual of
$\mbox{H}_*^{t+1} (\PP^N, \ii_X)$.
In fact, by the local duality
theorem for graded modules (see e.g. \cite[Theorem 3.6.19]{BH}),
the canonical module $\Omega_X$ is the graded
$k$-dual of $\mbox{H}_{\mathfrak m}^{t+1}(R_X)$, where $\mbox{H}_{\mathfrak m}^i$ denotes the $i$-th
local cohomology group with respect to the irrelevant maximal ideal $\mathfrak m$ of $R$,
and moreover
$\mbox{H}_{\mathfrak m}^{t+1}(R_X)\cong \mbox{H}^{t+1}_*(\mathbb{P}^N,\ii_X)$ (see
\cite[p.137]{HA}).
\end{proof}
For AG subschemes we have more ways to compute the invariant $m$:
\bp \label{bb1}
  Let $Z$ be an AG subscheme of $\PP^N$. Then $m(Z)$ is the integer $m$ such that $\Omega_Z \cong R_Z (m)$.
  Furthermore, if the minimal free resolution of $R_Z$ over $R$ has the form
  $$
  0 \ra R(-c) \ra \cdots \ra R \ra R_Z \ra 0,
  $$
  then $m(Z) = c-N-1$ and $b(Z)=c- \mbox{\em codim}(Z)$.
\ep
\begin{proof}
It is clear that $\Omega_Z \cong R_Z (m)$ implies
$m=m(Z)$. To relate the integer $c$ appearing at the last step of the resolution
with $m$, we compute
$\Omega_Z \cong \mbox{Ext}^{\mbox{\small codim $(Z)$}}(R_Z,\Omega_{{\mathbb P}^N})$
using the minimal free resolution of $R_Z$, and we find $$\Omega_Z \cong R_Z(c-N-1).$$
Hence $m = c-N-1$.
\end{proof}
\br
For an AG subscheme $Z$, the integer $b(Z)$ is called the {\em socle degree} of $R_Z$ in Migliore's book
\cite{M} and is denoted by $r$; the integer $m(Z)$ is sometimes referred to as the $a$-invariant
of the Gorenstein graded algebra $R_Z$.
\er

\bc\label{lemma32}
Let $X$ be a codimension 2 ACM  subscheme of ${\mathbb P}^N$.
Then in the minimal free resolution of $I_X$:
$$
0\longrightarrow \oplus R(-b_j) \longrightarrow  \oplus R(-a_i)  \longrightarrow
I_X \longrightarrow 0
$$
we have $\max \{b_j\}=b(X)+2$.
\ec
\begin{proof}
The regularity of $\ii_X$ equals $\max \{b_j\}-1$ by \cite[p. 8]{M} and $b(X)+1$ by Proposition \ref{bb}.
\end{proof}

Macaulay has characterized the possible
$h$-vectors for ACM subschemes of $\PP^N$ \cite{Macaulay}. We will only need this result for subschemes
of codimension 2 (Proposition \ref{C2}) and its analogue for AG subschemes of codimension 3
(Proposition \ref{G3}), which is due to Stanley.
\begin{defn} \label{C2def}
An $h$-vector is said to be $C2$-admissible if
there exists an $s \geq 1$ such that
$$
\begin{cases} h(n) = n+1 \;\;\;\mbox{if \ } 0 \leq n \leq s-1, \\
h(n) \geq h(n+1) \;\;\;\mbox{if \ } n \geq s-1. \end{cases}
$$
A $C2$-admissible $h$-vector is said to be  of decreasing type if
$h(a)>h(a+1)$ implies that for each $n \geq a$ either $h(n)>h(n+1)$ or $h(n)=0$.
\end{defn}

\bp \label{C2} \
\begin{enumerate}
\item
A finitely supported numerical function $h: \N \ra \N$ is the
$h$-{vector} of a codimension 2 ACM subscheme  $X$ of $\PP^N$ ($N \geq 2$) if and only if
$h$ is $C2$-admissible. Furthermore, $X$ can be taken
reduced, and a locally complete intersection in codimension $\leq 2$.
\item
If $X$ is an integral codimension 2 ACM subscheme of $\PP^N$, then $h_X$ is of decreasing type. Conversely,
if $h$ is of decreasing type and $N \geq 3$, there exists an integral codimension 2 ACM subscheme  $X \subset \PP^N$
with $h_X=h$ which is  smooth in codimension $\leq 2$ (thus smooth if $N=3,4$).
\end{enumerate}
\ep
\begin{proof}
The first statement is a very special case
of a theorem by Macaulay (see \cite{Macaulay},\cite{S}), and is equivalent, in case  $N=3$, to Theorem V.1.3 of \cite{MDP},
and, in general, to Proposition 1.3 of \cite{N}. The fact that $X$ can be taken
reduced, and a locally complete intersection in codimension $\leq 2$ is proven in \cite[3.2]{N}.
The second statement is proven in  \cite[3.3]{N};
it was first proven over the complex numbers and for $N=3$ in \cite[2.2 and 2.5]{GP}.
\end{proof}

\bigskip
We have also the result of Ellingsrud \cite{E} about the Hilbert scheme:

\noindent \bt[Ellingsrud]  The set of all ACM codimension $2$
subschemes of ${\mathbb P}^N$ ($N\geq 3$) with a given $h$-vector is a smooth,
open, irreducible subset of the Hilbert scheme of all
closed subschemes of ${\mathbb P}^N$. (There is also an explicit
formula for its dimension). \et

Stanley \cite{S}, drawing on Macaulay's Theorem, and applying the
structure theorem of Buchsbaum and Eisenbud \cite{BE},
characterized the possible $h$-vectors of AG codimension $3$
subschemes. Before we state his result, we define the ``first
half" $k_h$ of an $h$-vector $\{1,h(1), \ldots, h(b)\}$ setting:
\begin{equation}\label{key}
 k_h(n) = \left\{ \begin{array}{cl}
\partial h(n)  &\mbox{\rm for $0 \leq n \leq [b/2]$,} \\
0 &\mbox{\rm otherwise,}
\end{array} \right.
\end{equation}

where $[b/2]$ denotes the integral part of $b/2$. We have (see \cite[4.2]{S})

\bp  \label{G3}
A finitely supported numerical function $h: \N \ra \N$
is the $h$-vector of an AG codimension $3$ subscheme of ${\mathbb
P}^N$ ($N \geq 3$) if and only if
\begin{enumerate}
\item $h$ is symmetric, meaning
that $h(n) = h(b-n)$ for all $ 0 \leq n \leq b$;
\item the first half of $h$ is $C2$-admissible.
\end{enumerate}
We say such an $h$-vector is $G3$-admissible.
\ep

\br\label{remark23}
If $Z$ is a codimension 3 AG-subscheme of $\PP^N$, we let $k_Z$ denote
the first half of $h_Z$ and call it the $k$-vector of $Z$. Note that
the degree of $k_Z$ is precisely $h_Z ([b(Z)/2])$. Furthermore,
because of its symmetry, the $h$-vector of $Z$ is determined by its
first half $k_Z$ and by $b(Z)$.
\er

\bex
The $h$-vectors
$$
\{1,3,3,\ldots, 3,3, 1\}
$$
all are $G3$-admissible with first half $\{1,2\}$ - no matter how many $3$'s appear in
the string.
\eex

\br
There is an analogue of Ellingsrud's Theorem in the case of
codimension 3 AG-subschemes. If we fix the $h$-vector, the set
$\mbox{PGor}(h)$ of codimension 3 AG-subschemes of $\mathbb{P}^N$
with the given $h$-vector carries a natural scheme structure,
which makes it into a smooth irreducible subscheme of the Hilbert
scheme; furthermore, the dimension of $\mbox{PGor}(h)$ can be
computed, see \cite{D},\cite{K},\cite{KM}. Corollary \ref{dimension} below allows one
to compute $\dim  \mbox{PGor}(h)$ inductively as a function of the $h$-vector.
\er

We will need formulas for the variation of the $h$-vector under
liaison and biliaison. Recall the definition of Gorenstein liaison
\cite[5.1.2]{M} and biliaison \cite{H1}:

\bd Let $V_1$, $V_2$, $X$
be equidimensional subschemes of $\PP^N$, all of the same
dimension, with $V_1$ and $V_2$ contained in $X$. We say that
$V_2$ is {\em $G$-linked} to $V_1$ by $X$ if $X$ is AG and
$\ideal{V_2,X} \cong \Hom(\coo_{V_1},\coo_{X})$; in this case, it
is also true that $V_1$ is $G$-linked to $V_2$ by $X$.
\ed

\bd\label{defbil}
Let $V_1$ and $V_2$ be equidimensional
closed subschemes of dimension $t$ of ${\mathbb P}^N$.  We say
that $V_2$ is obtained by an {\em elementary  biliaison} of height
$h$ from $V_1$ if there exists an ACM scheme $X$ in ${\mathbb
P}^N$, of dimension $t+1$, containing $V_1$ and $V_2$ so that
$\ideal{V_2,X} \cong \ideal{V_1,X}(-h)$. (In the language of generalized divisors
on $X$, this says that
$V_2 \sim V_1 + hH$, where $H$
denotes the hyperplane class).
\ed

If $h\geq 0$ (respectively $h\leq 0$), we will speak of an ascending
(respectively descending) biliaison.
If we restrict the scheme $X$ in the definition to be
a complete intersection scheme, we will speak of {\em
CI-biliaison}.

The following proposition is well known and shows how the $h$-vector of an ACM scheme
changes under $G$-liaison and elementary biliaison.

\bp \label{variation} \

\begin{enumerate}
    \item Suppose $V_1$ and $V_2$ are two ACM subschemes of $\PP^N$ linked by the AG scheme $X$.
Let $b=b(X)$. Then
\begin{equation*}
h_{V_2}(n)=h_{X}(n)-h_{V_1}(b-n) \mbox{\ \ \ \ for all \ \ }n\in \Z.
\end{equation*}

    \item

Suppose $V_1$ and $V_2$ are ACM subschemes of $\PP^N$  such that
$V_2$ is obtained by an elementary  biliaison of height
$1$ from $V_1$ on the ACM scheme $X$. Then
\begin{equation*}
h_{V_2}(n)=h_{X}(n) + h_{V_1}(n-1) \mbox{\ \ \ \ for all \ \ }n\in \Z.
\end{equation*}

\end{enumerate}
\ep
\begin{proof}
The first statement is proven in \cite[5.2.19]{M}; the second follows from the isomorphism
$\ideal{V_2,X} \cong \ideal{V_1,X}(-1)$.
\end{proof}

\bex
We will later need the $h$-vector of a zero dimensional subscheme $W \subset \PP^3$ consisting of $d$ general points.
Since points in $W$ impose independent conditions on surfaces of degree $n$, the least degree
$s=s(W)$ of a surface containing $W$ is the unique positive integer such that $\binom{s+2}{3} \leq d < \binom{s+3}{3}$,
and
we have
\begin{equation}\label{hgeneral}
 h_W (n) =
\begin{cases}
\begin{array}{ll}
\ds \binom{n+2}{2}  & \mbox{\ \ if $0 \leq n<s$} \\
d- \ds \binom{s+2}{3}  & \mbox{\ \ if $n=s$} \\
0 & \mbox{\ \ otherwise.}
\end{array}
\end{cases}
\end{equation}
\eex

\section{AG divisors on codimension 2 ACM subschemes}
\label{sec3}

If $Z$ is an AG subscheme of codimension 3 in ${\mathbb P}^N$,
we know that its ``first half'' $k_Z$ is the $h$-vector
of a codimension 2 ACM subscheme $X \subset \PP^N$. Thus given a pair of $h$-vectors
$h$ and $k$, with $k$ equal to the first half of $h$, it is natural to
ask whether we can find a pair of subschemes $Z$ and $X$ of $\PP^N$ with $h_Z=h$ and $h_X=k$
satisfying further $Z \subset X$. Theorem \ref{thm25} gives a positive answer.

Consider an ACM scheme $X \subset \PP^N$ which is generically Gorenstein, so that in the language of
generalized divisors~\cite{H1} one can speak of the anticanonical
divisor $-K$ on $X$.
Then it is well known that any divisor $Z$ on $X$
linearly equivalent to $-K+mH$, that is,
such that $\ideal{Z,X} \cong \omega_X (-m)$,
is AG (\cite[5.2]{KMMNP}, \cite[3.4]{H1}).
In the following proposition we recall this fact, and describe the relation between the $h$-vectors of
$X$ and $Z$.

\bp \label{3.4}
Let $X$ be a generically Gorenstein ACM  subscheme ${\mathbb P}^N$ of dimension $t \geq
1$, and let $H$ be a hyperplane section of $X$.
Suppose $Z \subset X$  is  divisor linearly equivalent to $-K+ m H$,
for some integer $m$. Then
\begin{enumerate}
\item[(a)]
$Z$ is an AG subscheme of $\PP^N$ with $\Omega_Z \cong R_Z (m)$;
in particular, $m=m(Z)$.
\item[(b)] \label{parttwo}
 The $h$-vector of $Z$ is determined by the integer $m$ and the $h$-vector of $X$
by the formula
$$
\partial h_Z (n)= h_X (n)- h_X(m+t+1-n).
$$
In particular, $$b(Z) \geq b(X)$$
and, if $k_Z$ denotes the first half of $h_Z$, we have
$h_X (n) \geq k_Z(n)$ for every $n$.
\item[(c)]
The equality $h_X=k_Z$ holds if and only if $ b(Z) \geq 2 b(X)$
(or equivalently  $ m \geq 2 b(X) - t$),
so that in this case the $h$-vector of $X$ is determined by that of $Z$.
\end{enumerate}
\ep
\br\label{finite} \
  Since $b(Z) \geq b(X)$, we see that if $Z$ is fixed, there are at most finitely many
  $h$-vectors $\tilde{h}$ with the property that there exists an ACM $X$ with $h_X=\tilde{h}$
  on which $Z$ is a divisor linearly equivalent to $-K_X+ m H_X$.
\er

\noindent
{\em Proof of 3.1.}
{{From}} the exact sequence
\begin{equation}\label{antican}
\exact{ \Omega_X (-m)}{R_X}{R_Z}
\end{equation}
it follows
\newcommand{\depth}{\mbox{depth\ }}
$$
\depth R_Z \geq \mbox{Min} (\depth \Omega_X - 1, \depth R_X) = \dim R_X -1.
$$
Since $Z$ is a divisor on $X$, $\dim Z = \dim X -1$, so we must have equality, and thus $Z$ is ACM.

To see $Z$ is in fact $AG$, we apply the functor $\mbox{Hom}_{R_X} ( - , \Omega_X)$ to the exact
sequence (\ref{antican}) to get
$$
0 = \mbox{Hom}_{R_X} ( R_Z , \Omega_X) \ra \Omega_X \cong I_{Z,X} (m) \ra R_X (m) \ra \Omega_Z
\ra 0.
$$
This implies $\Omega_Z \cong R_Z (m)$ because $\Omega_Z$ is a faithful $R_Z$-module.

To prove (b), we use the exact sequence (\ref{antican}) and Serre duality to obtain
\begin{equation} \label{half}
\partial h_Z (n) = h_X(n) - \partial^{t+1} (h^t (\coo_X (m-n)).
\end{equation}
To see $\partial^{t+1} (h^t (\coo_X (m-n))=h_X(m+t+1-n)$, we use the fact
that $\partial^{t+1} h^t \coo_X (n)=(-1)^{t+1}\partial^{t+1}  h^0 \coo_X (n) $ because
$\chi \coo_X (n)= h^0 \coo_X(n) + (-1)^t h^t \coo_X (n)$ is a polynomial of degree
$t$ in $n$, together with the fact that for any numerical function $f(n)$,
if $g(n)=f(m-n)$, then
$\partial^{t+1} g(n)= (-1)^{t+1} \partial^{t+1} f(m+t+1-n)$.

Setting $n=0$ in the formula, we see that
$h_X(m+t+1)=0$, hence $ b(X)<m+t+1$. But $m+t=b(Z)$ because $m=m(Z)$ and $\dim Z=t-1$.
Thus $b(X) \leq b(Z)$.

As for part (c), it is clear that $h_X=k_Z$ implies $b(Z) \geq 2b(X)$. Conversely,
suppose  $b(Z) \geq 2 b(X)$.
Then $[b(Z)/2] \geq b(X)$, which implies $h_X (n)= 0 = k_Z (n)$ for $n > [b(Z)/2]$.
On the other hand, for $n \leq [b(Z)/2]$,
we have $m+t+1-n=b(Z)+1-n> b(X)$, hence by part (b)
$$k_Z (n)= \partial h_Z (n)= h_X (n).$$
\rightline{$\Box$}

\begin{prop}\label{mainthe}
Let $X$ be a codimension 2 ACM  subscheme of ${\mathbb P}^N$,
$N\geq 3$. Assume $X$ is a locally complete intersection in
codimension $\leq 2$. Then
\begin{enumerate}
\item[(a)] if $m\geq 2b(X)-N+1$,  there is an effective generalized divisor $Z \subset X$
linearly equivalent to $mH-K$;
\item[(b)] if $m\geq 2b(X)-N+2$ and furthermore $X$ is a
locally complete intersection everywhere,  the invertible sheaf
$\mathcal{L}(mH-K)= \omega_X^\vee(m)$ is very ample.
\end{enumerate}
\end{prop}
\begin{proof}
We will follow the proof of Theorem 3.5 in \cite{hartshorne1}.
First of all note that the assumption $X$ is a locally complete
intersection in codimension $\leq 2$ assures that $X$ satisfies
the conditions "$G_1+S_2$" under which the theory of generalized
divisors is developed in \cite{GD}; also $K_X$ is an almost
Cartier divisor on $X$, and $\mathcal{L}(mH-K)= \omega_X^\vee(m)$.

Let us begin by proving part (a): since $\ii_X (b+1)$ is generated
by global sections, so are $(\ii/\ii^2)(b+1)$ and
$\bigwedge^2(\ii/\ii^2) \,(2b+2)$. By assumption, off a closed
subscheme $Y$ of $X$ of codimension at least $3$, $X$ is a locally
complete intersection, so $\bigwedge^2(\ii/\ii^2)\cong
\omega_X^\vee(-N-1)$ on $X \setminus Y$, and
$\omega_X^\vee(2b-N+1)$ is generated by its global sections on $X
\setminus Y$. In particular, for $m \geq 2b-N+1$ we can find a
section of $\omega_X^\vee(m)$ over $X \setminus Y$ which generates
$\omega_X^\vee(m)$ at the generic points of $X$. Since $\codim \,Y >1$, we can
extend this section to all of $X$ \cite[1.11]{GD} . We conclude that for $m\geq
2b-N+1$ there is an effective divisor $Z \subset X$
linearly equivalent to $mH-K$ \cite[2.9]{GD}.

If $X$ is a locally complete intersection everywhere, then the
above argument shows that for $m\geq 2b-N+1$ the invertible sheaf
$\omega_X^\vee(m)$ is globally generated. Statement (b) now
follows from this and the
fact that the tensor product of a very ample invertible sheaf and
a globally generated invertible sheaf is very ample.
\end{proof}

\begin{thm} \label{thm25}Let $h$ be a $G3$-admissible $h$-vector,
and let $k$ be its first half as defined in section 2 formula (\ref{key}).
Denote by $b$ and $\overline{b}$ the largest integers such that $h(n)\not=0$ and
$k(n)\not=0$ respectively, and let $m=b-N+2$, $N\geq 3$. Then:
\begin{enumerate}
\item[(a)]
there is an AG subscheme $Z$ of $\mathbb{P}^N$ of codimension 3
with $h$-vector $h_Z=h$ lying on an ACM  scheme $X$ of codimension
2 with $h$-vector $h_X = k$;
\item[(b)]
if the first half $k$ of $h$ is of decreasing type, we may
take $X$ to be integral, and in case $N=3$ or $4$, both $X$ and $Z$ to be smooth;
\item[(c)]
if $b \geq 2 \overline{b}+2$ (that is, $m \ge
2\overline{b}-N+4$), then there is an open subset $V$ of the Hilbert scheme
$\mbox{PGor}(h)$ such that every $Z \in V$ is of the form $Z \sim
mH-K$ on a codimension 2 ACM scheme $X$ with $h$-vector $k$.
\end{enumerate}
\end{thm}
\br
Note that in the statement the integers $b$ and $\overline{b}$ depend only on $h$,
and by definition of the ``first half'' of $h$ we always have $b \geq 2 \overline{b}$.
Thus the hypothesis in $(c)$ holds unless $b= 2 \overline{b}$ or $b= 2 \overline{b}+1$.
\er

\br
In \cite{DNV} it is shown that an integral AG subscheme of codimension $3$ of $\PP^N$
exists with given $h$-vector if and only if the first half of that $h$-vector is of decreasing type.
A result related to statement $(c)$ above is Theorem 3.4 of \cite{B}; see also Remark 5.3 of \cite{KMMNP}.
For a deformation theoretic approach to the question of determining those AG schemes $Z$ of the form
$Z \sim mH_X-K_X$ see \cite{Kl}.
\er

\noindent
\begin{proof} (a) Suppose given a $G3$-admissible $h$-vector $h$. Then its first half $k$
is $C2$-admissible, hence there exists \cite[3.2]{N} an ACM subscheme $X \subset \PP^N$ of codimension 2
with $h$-vector $h_X=k$ which is
reduced and a locally complete
intersection in codimension $\leq 2$.
Then $b(X)=\overline{b}$, and by definition of
the ``first half'' function $b \geq 2 \overline{b}$, that is,
$m \ge 2\overline{b}-N+2$. By Proposition \ref{mainthe}, there is a divisor $Z \subset X$ linearly equivalent
to $mH-K$, and, by Proposition \ref{3.4}, $Z$ is AG with $m(Z)=m$, hence $b(Z)=b$ and the first half
of  the $h$-vector $h_Z$ equals $h_X$. Thus $h_Z$ and $h$ both have first half $k$ and are last nonzero
at $b$, hence by symmetry $h_Z=h$.

\noindent
(b) If $k$ is of decreasing type, then $X$ can be taken to be integral \cite[3.3]{N}.
If $N=3$ or $4$, then $X$ can be taken even smooth \cite[3.3]{N}, and,
since $mH-K$ is very ample by Proposition \ref{mainthe}, we can take $Z$ to be smooth by the usual Bertini
theorem.
\noindent

\noindent
(c) It is  known that
$\mbox{PGor}(h)$ is smooth, and the tangent space to
$\mbox{PGor}(h)$ at the point corresponding to the subscheme $Z$
is isomorphic to the degree zero piece $\ _0 \rmhom_R(I_Z,R_Z)$ of
the graded $R$-module $\rmhom_R(I_Z,R_Z)$: see for example
\cite{KM}, and \cite{K} for the zero-dimensional case. We remark
that, in case $\dim Z > 0$, the tangent space $\ _0
\rmhom_R(I_Z,R_Z)$ is isomorphic to the space of global sections
of the normal sheaf $\cnn_Z$ of $Z$ in $\PP^N$.

It is therefore sufficient to show that the dimension of
the family  in statement (c) is greater or equal than the
dimension of $\ _0 \rmhom_R(I_Z,R_Z)$.

We will restrict $X$ to vary in the open subset of locally  complete intersections
in codimension $\leq 2$. Therefore we can use the theory of generalized
divisors on $X$. Since $X$ is an ACM scheme and since the linear system
$|Z|=|mH-K|$ is effective,
the dimension of the linear system
$|Z|$ on $X$ is equal to $\dim \zhomx {I_{Z,X}}{R_X}-1$ \cite[2.9]{GD}.

Since $I_{Z,X} \cong \Omega_X (-m)$, we have
$\rmhom_{R_X}(I_{Z,X},I_{Z,X}) \cong R_X$ and
$\ext^1_{R_X}(I_{Z,X},I_{Z,X})=0$. Thus applying the functor
$\rmhom_{R_X}(I_{Z,X}, - )$ to the sequence
\begin{equation}\label{seqJ}
0 \lra I_{Z,X} \lra R_X \lra R_Z \lra 0
\end{equation}
we obtain an exact sequence
$$
0 \lra R_X \lra \rmhom_{R_X}(I_{Z,X},R_X)\lra
\rmhom_{R_X}(I_{Z,X},R_Z)\lra 0.
$$
Hence
\begin{equation} \label{dls}
\dim \ _0\rmhom_{R_X}(I_{Z,X},R_Z)   =\dim \ _0
\rmhom_{R_X}(I_{Z,X},R_X)-1
\end{equation}
is the dimension of the linear system $|Z|$ on $X$.

The family of ACM schemes in which $X$ varies has dimension
\begin{equation} \label{dx}
\dim \zhomr{I_X}{R_X}
\end{equation}
because codimension 2 ACM subschemes of positive dimension are
unobstructed~\cite{E}.

Next we remark that if $Z$ is as in $(c)$ then by Proposition \ref{3.4}
the $h$-vector $h_X$ equals $k_Z$. The ACM scheme $X$ is unique. Indeed,
the equality of the Hilbert functions of $X$ and $Z$ for $n\leq \overline{b}$
guarantees that the map $H^0(\ii_X (n)) \to H^0(\ii_Z (n))$ is an
isomorphism for $n\leq \overline{b}$. Since $\ii_X$ is generated in degrees
less or equal than $\overline{b}+1$, the homogeneous ideal of $X$ is uniquely
determined by the ideal of $Z$.

Therefore the dimension of the family in statement $(c)$ is the
sum of the dimension of the linear system $|Z|$ on $X$ and of
dimension of the Hilbert scheme in which $X$ varies.
By~(\ref{dls}) and~(\ref{dx}) this dimension is equal to
$$
\dim \zhomr{I_X}{R_X}  + \dim \ _0\rmhom_{R_X}(I_{Z,X},R_Z) .
$$
Thus to complete the proof it suffices to show
\begin{equation}\label{df}
\dim \ _0 \rmhom_{R} (I_Z,R_Z) \leq \dim \zhomr{I_X}{R_X}  + \dim
\ _0\rmhom_{R_X}(I_{Z,X},R_Z).
\end{equation}

To this end, we apply $\rmhom_R(-,R_Z)$ to the exact sequence
\begin{equation} \label{seq2}
0  \lra I_X \lra I_Z \lra I_{Z,X} \lra 0
\end{equation}
and obtain a new exact sequence
$$
0 \lra \homr {I_{Z,X}}{R_Z} \lra \homr {I_Z}{R_Z} \lra \homr
{I_X}{R_Z}, $$ from which we deduce:
\begin{equation}\label{dg}
\dim \ _0\homr {I_Z}{R_Z} \leq \dim \ _0\homr {I_X}{R_Z}   + \dim
\zhomr {I_{Z,X}}{R_Z} .
\end{equation}

Comparing~(\ref{df}) and (\ref{dg}), we see we will be done if we
can show $ \homr {I_X}{R_X}  \cong \homr {I_X}{R_Z}$. For this, we
apply the $\rmhom_R(I_X,-)$ to the exact sequence (\ref{seqJ}) and
obtain
$$
0\lra \homr {I_X} {I_{Z,X}} \lra \homr {I_X} {R_X} \lra \homr
{I_X} {R_Z} \lra  \ext^1_R( I_X, I_{Z,X}).
$$
Since $I_{Z,X} \cong \Omega_X(-m))$, we need to
show that $ \homr {I_X} {\Omega_X(-m)} = \ext^1_R( I_X,
\Omega_X(-m)) =0$.

To prove these vanishings, we consider the minimal presentation of
$I_X$ over $R$:
\begin{equation}\label{seqI2}
\bigoplus R(-b_j) \lra \bigoplus R(-a_i) \lra I_X \lra 0.
\end{equation}

The module $\rmhom_{R}(I_X, \Omega_X (-m))$ is a submodule of
$\bigoplus \Omega_X (a_i-m)$. By Proposition~\ref{lemma32} we have
$\overline{b}+2 = \mbox{Max} \{b_j\} \geq a_i $ for every i. Since by
assumption $m \geq 2\overline{b}-N+4$, we have
$$
m-a_i \geq \overline{b} -N +2 = \overline{b} + \dim X.
$$
By Proposition~\ref{bb} it follows $\Omega_X (a_i-m)= 0$ for
every $i$, hence $\rmhom_{R}(I_X, \Omega_X (-m))=0$.

To deal with $\ext^1_{R}(I_X,\Omega_X(-m))$, let $M$ denote the
kernel of the surjection $\bigoplus R(-a_i) \ra I_X$. Then
$\ext^1_{R}(I_X,\Omega_X(-m))$ is a quotient of $\rmhom_{R}(M,
\Omega_X (-m))$, which in turn is a submodule of $\bigoplus
\Omega_X (b_j-m))$. Since $m \geq 2\overline{b}-N+4$, the same argument as
above shows $\Omega_X (b_j-m)=0$ for every $j$. Hence
$\rmhom_{R}(M, \Omega_X (-m))=0$, and, a fortiori,
$\ext^1_{R}(I_X,\Omega_X(-m))=0$, which is what was left to prove.
\end{proof}

\section{Arithmetically Gorenstein zero-dimensional subschemes}
A natural problem is to look for geometric conditions which
allow one to decide whether or not a subscheme of $\pp^N$ is AG. In
the case of zero-dimensional subschemes there is a characterization
in terms of the so called Cayley--Bacharach property:
it is the implication $(a) \Longleftrightarrow  (b)$ in the following proposition,
and it was originally proven in \cite{DGO}
in the reduced case and then generalized to the non reduced
case in \cite{Kr}. Below we give an alternative proof.

Recall a zero dimensional subscheme $Z \subset \PP^N$ is ACM, and
by Proposition \ref{bb} its $m$-invariant is
$$m(Z)=\max\{ n\ |\ (\Omega_Z)_{-n}\not=0\} = \max\{ n\ :\  h^{1} (\PP^N, \ii_Z (n) )\not=0\},$$
and the regularity of $Z$ is $m(Z)+2$.
\bp \label{CB}
Let $Z$ be a zero-dimensional scheme in $\pp^N$, let $m=m(Z)$ and suppose the $h$--vector of $Z$
is symmetric, i.e.
$$h(n)=h(m+1-n) \mbox{\ \ \ \ \ for every $n \in \Z$}.$$
Then the following are equivalent:
\begin{enumerate}
\item[(a)] $Z$ is AG;
\item[(b)] $Z$ satisfies the Cayley--Bacharach property:
for every subscheme $W\subset Z$ with $\deg (W)=\deg (Z)-1$,
we have
$$h^0(\pp^N,\mathcal{I}_W (m))=h^0(\pp^N,\mathcal{I}_Z(m));$$
\item[(c)] for any proper subscheme $W_0$ of $Z$, $m(W_0)<m$;
\item[(d)] for any proper subscheme $W_0$ of $Z$,
$h^{1} (\PP^N, \ii_{W_0} (m) )=0$.
\end{enumerate}
\ep
\begin{proof} \

\vspace{2mm}
(c)$\Longleftrightarrow$ (d)
Since $ \bigoplus_{n \in \Z} \mbox{H}^{1} (\PP^N, \ii_{W_0} (n) )$ is the graded
$k$-dual of $\Omega_{W_0}$, the vanishing $$h^{1} (\PP^N, \ii_{W_0} (m) )=0$$ implies
$h^{1} (\PP^N, \ii_{W_0} (n) )=0$ for every $n \geq m$, that is, $m (W_0) < m$.

\vspace{2mm}
(b)$\Longleftrightarrow$ (d)
Since $m(W_0) \leq m(W)$ whenever $W_0 \subseteq W \subseteq Z $, we may assume in (d)
that $W_0=W$ has length  $\deg Z-1$.

Looking at Hilbert polynomials we see
$$
h^0 (\ii_W (m)) - h^1 (\ii_W (m)) = h^0 (\ii_Z (m)) - h^1 (\ii_Z (m))+1.
$$

The assumption on the $h$-vector implies  that
$ h^1 (\PP^N, \ii_Z (m)) = h_Z (m+1) =1$. Thus
$ h^0 (\ii_W (m))  = h^0 (\ii_Z (m))$ if and only if $h^1 (\ii_W (m))=0$.

\vspace{2mm}
(c)$\Longrightarrow$ (a)
The symmetry of the $h$-vector implies that $R_Z$ and $\Omega_Z(-m)$
have the same Hilbert function.

We can define a function $\varphi:\ R_Z\to \Omega_Z(-m)$ sending
$1$ to $\alpha$ where $\alpha$ is a non zero element of degree zero in $\Omega_Z(-m)$.
If $\varphi$ had a
nontrivial kernel,
then it would define a proper subscheme $W$ of $Z$ together with
a non zero map $R_W \to \Omega_Z(-m)$ and so $m(W)=m(Z)$,
contradicting (c). Thus $\varphi$ has trivial kernel, hence it is an isomorphism
because $R_Z$ and $\Omega_Z(-m)$ have the same Hilbert function. Therefore $Z$ is AG.

\vspace{2mm}
(a)$\Longrightarrow$ (c)
Suppose $m(W)=m$. Then we have a nontrivial map
$$
R_W \ra \Omega_W(-m) \ra \Omega_Z(-m) \cong R_Z
$$
lifting the natural map $R_Z \ra R_W$, hence $W=Z$.
\end{proof}

For later use we need to rephrase this criterion in the case $Z$ is a
subscheme of an AG curve.
Thus suppose we are given an AG curve $C$ and a zero dimensional
subscheme $Z \subset C$. Then $Z$ can be thought of as an effective
generalized divisor on $C$ \cite{GD}, and it corresponds to a nondegenerate
section of the reflexive sheaf of $\coo_C$ modules
$$
\coo_{C}(Z) = \mathcal{H}om_{\coo_C} (\ideal{Z,C}, \coo_C).
$$
Analogously we define
$$
R_C (Z) = \mbox{Hom}_{R_C} (I_{Z,C}, R_C).
$$
Note that, since $C$ is ACM, we have
$I_{Z,C} \cong \mbox{H}^0_* (\ideal{Z,C})$ and
$R_C \cong \mbox{H}^0_* (\coo_C)$, hence the $n$-th graded piece
$R_C (Z)_n$ of $R_C (Z)$ is isomorphic to $\mbox{H}^0
(\coo_C(Z+nH))$, which at the level of graded modules we may
rewrite as
$$
R_C (Z) \cong \mbox{H}^0_* (\coo_C(Z)).
$$

Applying the functor  $\mbox{Hom}_{R_C} ( - , R_C)$ to the exact
sequence
$$
\exact{I_{Z,C}}{R_C}{R_Z}
$$
we obtain another exact sequence, analogous to \cite[2.10]{GD}
\begin{equation}\label{ggr}
\exact{R_C}{R_C(Z)}{\mbox{Ext}^1_{R_C} ( R_Z, R_C) \cong \Omega_Z \otimes \Omega_C^{\vee}}
\end{equation}
where the last isomorphism can be derived as follows:
$$
\mbox{Ext}^1_{R_C} ( R_Z, R_C) \cong \mbox{Ext}^1_{R_C} ( R_Z, \Omega_C) \otimes \Omega_C^{\vee}
\cong \Omega_Z \otimes \Omega_C^{\vee}.
$$
\bc \label{curvyCB}
Suppose $Z$ is an effective divisor on the AG curve $C$. Let
$m=m(Z)$ and $e=m(C)$.  Suppose the $h$--vector of $Z$ is symmetric.
 Then $Z$ is AG if and only if
$$
\dim (R_C (W))_{e-m}  < \dim (R_C (Z))_{e-m}
$$
for every subscheme $W\subset Z$ with $\deg (W)=\deg (Z)-1$.
\ec
\begin{proof}
Pick a subscheme $W\subset Z$ with $\deg (W)=\deg (Z)-1$.
Since $\Omega_C \cong R_C (e)$, applying (\ref{ggr}) to $W$ and then to $Z$,
we see that $m(W) < m$ is
equivalent to
$$\dim (R_C (W))_{e-m}  < \dim (R_C (Z))_{e-m}.$$
\end{proof}

\section{Complete Intersection Biliaison}
In this section,  we will show
that a general AG set of points in $\PP^3$ is obtained by ascending complete
intersection biliaisons from
a point (see Definition \ref{defbil}).  We follow closely section $4$ of~\cite{hartshorne1}.

\bl \label{41}
Let $Z$ be a codimension 3 AG subscheme of  ${\mathbb P}^N$ with
$h$-vector $h$.  Let $s$ denote the minimum degree of a hypersurface containing $Z$,
$b= b(Z) = \max\{n \mid h(n) > 0\}$, and let $m$ be the integer for which
$\Omega_Z \cong R_Z(m)$.  Then
\begin{itemize}
\item[a)]   $s = \min\{n > 0 \mid h (n) < \binom{n+2}{2}\}$.
\item[b)]  $m = b - \dim (Z) -1$.

\item[c)] ${\mathcal I}_Z(b+2-s)$ is generated by global sections.

\end{itemize}
\el
\begin{proof}
Part a) follows from the definition of the $h$-vector.
Part b) we recall for memory from Proposition \ref{bb}.  For part c) we use the
theorem of Buchsbaum--Eisenbud \cite{BE} in the notation of \cite[\S 5,
pp.~62-63]{HTV}.  Let $R$ be the homogeneous coordinate ring of ${\mathbb
P}^N$.  Then the homogeneous ideal $I_Z$ of $Z$ has a resolution of the form
\[
0 \rightarrow R(-c) \rightarrow \oplus R(-b_i) \rightarrow \oplus R(-a_i)
\rightarrow I_Z \rightarrow 0
\]
with $i = 1,2,\dots,2r+1$ for some positive integer $r$.  Moreover, this
resolution is symmetric in the sense that if we order $a_1 \le a_2 \le
\dots \le a_{2r+1}$ and $b_1 \ge b_2 \ge \dots \ge b_{2r+1}$, then $b_i = c -
a_i$ for each $i$.  Furthermore, if we let $u_{ij} = b_i - a_j$ be the
associated degree matrix, then $u_{ij} > 0$ for $i+j = 2r+3$.

To relate this to the invariants $s$ and $b$ of the $h$-vector,
first note that the $a_i$ are the degrees of a minimum set of
generators of $I_Z$. Hence $a_1 = s$, which is the least degree of
a generator. By symmetry, $b_1 =c-s$. On the other hand, $b = c-3$
by Proposition \ref{bb1}. {From} the inequality $u_{2,2r+1} > 0$ we
find $b_2 > a_{2r+1} = \max\{a_i\}$.  But $b_1 = c - s \ge b_2$,
so we find $\max\{a_i\} < c-s$.  Hence $\max\{a_i\} \le b+2-s$,
and ${\mathcal I}_Z(b+2-s)$ is generated by global sections.
\end{proof}

\bt \label{thm-bil}
For any  $h$-vector $h$ corresponding to
an AG zero dimensional subscheme of ${\mathbb P}^3$ (as in (\ref{G3})),
there is a nonempty open
subset $V_{h}$ of the Hilbert scheme $\mbox{\em PGor}(h)$  such that any
$Z \in V_{h}$ can be obtained by strictly
ascending CI-biliaisons from a point in ${\mathbb P}^3$.
\et

\begin{proof}
We will prove, by induction on the degree, the following slightly
more precise statement.  For each $h$, let $s=s(h)$ and $t= t(h)$
denote the value of the invariants $s(Z)$ and $m(Z)+3-s(Z)$
for $Z \in \mbox{PGor}(h)$. We claim there is an open set
$V_{h}\subseteq \mbox{PGor}(h)$ such that for any $Z \in V_{h}$
\begin{itemize}
\item[(i)] There is a reduced complete intersection curve $C = F_s \cap F_{t}$
such that $Z$ is contained in the smooth locus $C^{sm}$ of $C$ and intersects every
irreducible component of $C^{sm}$.
\item[(ii)] There is an AG zero-scheme $Z' \sim Z - H$ on $C$, with $h$-vector $h'$,
such that $Z' \in V_{h'}$.
\end{itemize}

To begin with, by Lemma \ref{41} a zero dimensional AG subscheme $Z \subseteq \PP^3$
is contained in the complete intersections $C$ of two surfaces of degree $s=s(Z)$ and
$t=m(Z)+3-s(Z)$ respectively. Thus  property (i) is an open condition
on $\mbox{PGor}(h)$.

We start the induction with AG subschemes $Z$ having $s = 1$.  These are contained
in a ${\mathbb P}^2$, so they are complete intersections, and for these
the theorem is well known.

So now we assume $s \ge 2$.  Suppose for a moment that $Z \subseteq C$ satisfies
condition (i).  We will show that the linear system $|Z-H|$ is nonempty and
contains an AG subscheme $Z'$.  We use the exact sequence (\ref{ggr})
twisted by $-H$:
\begin{equation}\label{ottobis}
0 \rightarrow R_C(-H) \rightarrow R_C(Z-H) \rightarrow
\Omega_Z \otimes \Omega_C^{\vee}(-H) \rightarrow 0.
\end{equation}
Now $\Omega_C \cong R_C (s+t-4) = R_C(m-1)$, and
$\Omega_Z \cong R_Z (m)$, so the sheaf on the right
is just $R_Z$.  Therefore
\begin{equation}\label{nove}
h^0(\coo_C(Z-H))= \dim (R_C(Z-H))_0 = \dim (R_Z)_0 =1
\end{equation}
so the sheaf $\coo_C (Z-H)$ has a unique section $\sigma$ whose
restriction to $Z$ is $1$. {From} the condition that $Z$ meets every
irreducible component of $C^{sm}$, and $C$ being reduced, we
conclude that $\sigma$ is nondegenerate, and defines an effective
divisor $Z' \sim Z - H$ \cite[2.9]{GD}.  Furthermore, since $\sigma$
restricted to $Z$ is $1$, we find that $\mbox{Supp}\, (Z) \cap
\mbox{Supp}\,(Z') = \emptyset$.

We claim that $Z'$ is AG. First of all from the equality
\begin{equation} \label{recipe}
h_{Z'} (n) = h_Z (n+1) - h_C(n+1)
\end{equation}
and the fact that both $h_Z$ and $h_C$ are $h$-vectors satisfying $h(n) = h(m+1-n)$
for all $n \in \Z$, we see that $h_{Z'}$ is symmetric and
$m(Z')=m-2$. Therefore to show $Z'$ is AG it is enough
by Corollary \ref{curvyCB}  to show that
for every $D \subset Z'$ with
$\deg (D) = \deg(Z')-1$ we have
$$\dim (R_C (D))_1 < \dim (R_C (Z'))_1.$$

Since $Z' \sim Z-H$ is a Cartier divisor on $C$, the divisor
$Z'-D$ is effective and has degree $1$, therefore it is a point $Q$
in the support of $Z'$. Now $D+H \sim Z-Q$ and $Z'+H \sim Z$, so
what we have to prove is that
$$h^0 (\coo_C (Z-Q))=\dim (R_C (Z-Q))_0 < \dim (R_C (Z))_0 =h^0 (\coo_C (Z)).$$
This follows from the fact that  $\mbox{Supp}\, (Z) \cap \mbox{Supp}\,(Z') = \emptyset$, hence the section of $\coo_C (Z)$ corresponding to $Z$ is not
the image of a section of $\coo_C (Z-Q)$.

Now we explain the induction step of the proof.  Given an admissible $h$-vector
$h$ with $s \ge 2$, define $h'$ as in (\ref{recipe}) above.  By the induction
hypothesis there exists an open set $V_{h'} \subseteq \mbox{PGor}(h')$ of
AG schemes satisfying (i) and (ii).  Let $Z'$ be such a scheme, and let $Z' \subseteq
C' = F_{s'} \cap F_{t'}$ satisfy (i).  Note that
we either have $s'=s-1$ and $t'=t-1$ or $s'=s$ and $t'=t-2$.
So define a curve $C = (F_{s'} + H_1) \cap (F_{t'} +
H_2)$ or $C = F_{s'} \cap (F_{t'} + H_1 + H_2)$, where $H_1$, $ H_2$ are
planes in general position.  Then $C$ is a reduced complete intersection
curve of two surfaces of degree $s$ and $t$ respectively.

On this curve $C$, we will show, by an argument analogous to the one above,
that a general divisor $Z$ in the linear system $Z'+ H$ on $C$ is AG.

We first prove that an effective divisor $Z \in |Z'+H|$ is AG if and only if
$\mbox{Supp}\, (Z) \cap \mbox{Supp}\,(Z') = \emptyset$.
By construction the $h$-vector of $Z$ is symmetric.
Therefore by Corollary \ref{curvyCB} the divisor $Z$ is AG if and only if for every
$W \subset Z$
of degree $\deg(Z)-1$ we
have
$$\dim (R_C (W))_{-1} < \dim (R_C (Z))_{-1}.
$$
Since $Z$ is Cartier on $C$, we may write $W=Z - Q$ where $Q$ is a
point in the support in $Z$, and then the inequality above is
equivalent to
$$
h^0 (\coo_C (Z'-Q)) < h^0 (\coo_C (Z')).
$$

Now from the exact sequence
\begin{equation}\label{gr}
\exact{R_C}{R_C(Z')}{\Omega_{Z'} \otimes \Omega^{\vee}_{C}\cong R_{Z'} (-1)}
\end{equation}
we deduce $h^0 \coo_C (Z')=1$. Thus the inequality above is satisfied
if and only if for every  $Q \in \mbox{Supp}(Z)$
we have $h^0 \coo_C (Z' -Q) =0$, that is, $Q$ is not in $\mbox{Supp}(Z')$.

This completes the proof of the claim that $Z$ is AG if and only if $\mbox{Supp}\, (Z) \cap \mbox{Supp}\,(Z') = \emptyset$.
Now we show that a general $Z$ in $|Z'+H|$ satisfies this property. For this, note
that $I_{Z,C} \cong I_{Z',C}(-1)$, hence twisting by one the exact
sequence (\ref{gr}) we obtain
\begin{equation}
\exact{R_C(1)}{R_C(Z)}{R_{Z'}}.
\end{equation}
The degree zero piece is
$$
\exact{\mbox{H}^0 \coo_C (H)}{\mbox{H}^0 \coo_C (Z)}{(R_{Z'})_0}
$$
which shows that a general section of $\mbox{H}^0 \coo_C (Z)$ maps to a unit in $(R_{Z'})_0$,
hence does not vanish at any point $P$ of $Z'$.

Furthermore, since the trivial biliaison $Z' + H$
satisfies (i), and this is an open condition, we can choose $Z$ in $|Z'+H|$ so that
it is AG and satisfies (i).

Thus there exists an open subset of AG subschemes $Z \in \mbox{PGor}(h)$ satisfying (i).
Since the procedures of constructing $Z'$ from $Z$ and $Z$ from $Z'$ are
reversible, we can find an open subset $V_{h} \subseteq \mbox{PGor}(h)$ of AG
schemes $Z$ satisfying (i) with the associated scheme $Z'$ lying in $V_{h'}$.

This completes the inductive proof of (i) and (ii).  To prove the theorem, we
take a $Z \in V_{h}$, and by (ii) find a $Z' \in V_{h'}$ with smaller
degree.  We continue this process until either the degree is $1$ or $s = 1$,
which we have discussed above.

\end{proof}
As a corollary we now derive a formula that allows one to compute the dimension of $\mbox{PGor}(h)$
inductively (cf. \cite{CV} for a different approach). Fix a $G3$-admissible $h$-vector $h$, and let
$s=s(h)=\min \{n>0: \, h (n) < \binom{n+2}{2}\}$ and $t=m(h)+3-s(h)$.
Denote by $h'$ the $h$-vector defined as in (\ref{recipe}) by the formula
$h'(n)=h(n+1)-h_{s,t}(n+1)$, where $h_{s,t}$ denotes the $h$-vector of a
a complete intersection of two surfaces of degree $s$ and $t$.
\bc\label{dimension}
With the above notation, if $s\geq 2$, then
$$\dim\mbox{\em PGor}(h)=\dim \mbox{\em PGor}(h')- h'(s)-h'(t)+st+3s+3-\varepsilon,$$
where $\varepsilon=0$ if $t>s$ and $\varepsilon=1$ if $t=s$.
\ec
\begin{proof}
As in the proof of Theorem \ref{thm-bil}, given  a general $Z \in \mbox{PGor}(h)$, there are
$Z' \in  \mbox{PGor}(h')$  and a curve $C$, complete intersection of surfaces of degree $s$ and $t$,
such that $Z$ is linearly
equivalent to $Z'+H$ on $C$. The $h$-vector $h_C=h_{s,t}$ of $C$ is given by  $h_C(n)=n+1$ for $0\leq n\leq s-1$,
$h_C(n)=s$ for $s-1\leq n\leq t-1$ and $h_C(n)=h_C(m+1-n)$ for $n\geq t$.
Let us denote by $\mathcal{F}$ the family of complete intersection
curves $C=F_s\cap F_t$ containing the AG zero-scheme $Z$
as a divisor $Z\sim Z'+H$ and by $\mathcal{F}'$ the family of complete intersection
curves $C=F_s\cap F_t$ containing the AG zero-scheme $Z'$. We have
\begin{equation}\label{dimPGor}
\dim\mbox{PGor}(h)=\dim\mbox{PGor}(h')+\dim \mathcal{F}'+\dim_C|Z'+H|-\dim \mathcal{F},
\end{equation}
because $\dim_C|Z'|=0$ by formula (\ref{nove}).
Let us begin by supposing $t>s$ and
let us compute $\dim \mathcal{F}$ and $\dim \mathcal{F}'$. We have
$$
\dim \mathcal{F}= (h^0 \mathcal{I}_Z(s)-1)+(h^0 \mathcal{I}_Z(t)-h^0 \mathcal{O}_{\mathbb{P}^3}(t-s)-1)=
$$
$$= \left({{s+3}\choose 3}- \sum_{n=0}^{s} h(n)-1\right)+
\left({{t+3}\choose 3}- \sum_{n=0}^{t} h(n)-{{t-s+3}\choose 3}-1\right).
$$
To determine $\dim \mathcal{F}'$ we need the Hilbert function
of $Z'$ which equals $\ds \sum_{r=1}^{n+1} (h(r)-h_C(r))$.
$$
\dim \mathcal{F}'= \left({{s+3}\choose 3}- \sum_{n=1}^{s+1} h(n)+\sum_{n=1}^{s+1}
h_C(n)-1\right)+
$$
$$
\left({{t+3}\choose 3}- \sum_{n=1}^{t+1} h(n)+\sum_{n=1}^{t+1} h_C(n)-
{{t-s+3}\choose 3}-1\right).
$$
The exact sequence (\ref{ottobis}) shows that $\dim_C |Z'+H|=4$.
By substituting the various pieces in (\ref{dimPGor}) we obtain:
$$
\dim\mbox{PGor}(h)=\dim\mbox{PGor}(h')-h(s+1)-h(t+1)+\sum_{n=0}^{s+1} h_C(n)+
\sum_{n=0}^{t+1} h_C(n)+4.
$$
Using the properties of the $h$-vectors $h_C$, $h'$ we obtain
$$
\dim\mbox{PGor}(h)=\dim\mbox{PGor}(h')-h(s+1)-h(t+1)+ h_C(s)+
 h_C(s+1)+st+3s+1=
$$
$$
=\dim\mbox{PGor}(h')-h'(s)-h'(t)+h_C(s)-
 h_C(t+1)+st+3s+1.
$$
Moreover $h_C(s)=s-1$ if $t=s$, $h_C(s)=s$ if $t\geq s+1$, so $h_C(s)-h_C(t+1)=1$ if $t=s$,
$h_C(s)-h_C(t+1)=2$ if $t\geq s+1$,
from which the statement follows.

When $t=s$, the calculation is similar, the only difference being that in this case
we have $\dim \mathcal{F}'=\dim \mbox{Grass}(2,H^0(\mathcal{I}_{Z'}(s)))$ and analogously for $\dim \mathcal{F}$.
\end{proof}
\br \label{pdim}
Using Corollary \ref{dimension} we can compute the dimension of $\mbox{PGor}(h)$ by induction
on $s$,
once we know the dimensions for all $h$ with $s=1$. Now a zero dimensional subscheme $Z$ with $s=1$
is the complete intersection of two plane curves of degree $p$ and $q$,
$q\geq p$,
and thus the dimension count is immediate:
$$
\dim \mbox{PGor}(h)=
\left\{\begin{array}{ll}
3,\ &{\it if}\ p=q=1\\
q+4,\ &{\it if}\ p=1, q>1\\
p^2+3p+1,\ &{\it if}\ q=p\geq 2\\
pq+3p+2,\ &{\it if}\ q>p\geq 2.\\
\end{array}
 \right.
$$
\er

\vspace{5mm}

Given a $C2$-admissible $h$-vector $h$ (see Definition \ref{C2def}),
an argument similar to the one in Corollary \ref{dimension} can be used to determine
the dimension of the family $\mbox{ACM}(h)$ of ACM curves in $\mathbb{P}^3$ with a fixed
$h$-vector $h$. This dimension was first computed by Ellingsrud \cite{E}.
Here we give a formula which allows to compute $\dim \mbox{ACM}(h)$ inductively
(compare \cite{MDP}, Proposition 6.8 p.176).
Let $s=s(h)$ denote the least degree of a
surface containing a curve $C$ in $\mbox{ACM}(h)$. Let $h'$ be the $h$-vector  defined by
\begin{equation}\label{hvectorbil}
h'(n)=\left\{
\begin{array}{ll}
h(n)\ \ &{\rm for}\ n\leq s-2 \\
h(n+1)\ \ &{\rm for}\ n\geq s-1.
\end{array}
\right.
\end{equation}
Note that $h'$ is the $h$-vector of an ACM curve $C'$ obtained from
$C\in\mbox{ACM}(h)$ by performing
an elementary descending biliaison of height one on a surface of degree $s$.
Then:

\bp\label{adim}
Let $h$, $h'$ as above. Then
\begin{equation}\label{dimACM}
\dim\mbox{\em ACM}(h)=\dim\mbox{\em ACM}(h')+4s+\sum_{n\geq s+1} h'(n).
\end{equation}
\ep

\bc\label{cordimACM}
If $h(n)=0$ for $n\geq s+2$, then $\dim\mbox{\em ACM}(h)=4d$ where $d=\sum_{n\geq 0} h(n)$
is the degree of a curve in $\mbox{\em ACM}(h)$.
\ec
\begin{proof} By induction on $s$, beginning with $s=1$, in which case
we have a plane curve of degree 1, 2, or 3 and the result is known.  Alternatively,
one can observe that the hypothesis on the $h$-vector implies that any curve
$C\in \mbox{ACM}(h)$ has index of speciality $m(C) \leq s-1$; hence
the normal bundle $N_C$ satisfies $h^1(N_C)=0$ by \cite[Lemma 4.2]{cW}, and this implies the statement
by deformation theory.
\end{proof}

\section{AG zero dimensional subschemes of  $\PP^3$}
In this section we study AG zero dimensional subschemes of  $\PP^3$ of low degree and their Hilbert schemes
$\mbox{PGor}(h)$. The results are summarized in Table 1.

We begin by addressing the question of how many general points one can impose
on a general $Z\in\mbox{PGor}(h)$.
We denote by $\mu=\mu(h)$ this number.

More generally, when $\caf$ is an irreducible flat family of subschemes of $\PP^N$,
we define $\mu(\caf)$ (or $\mu(Y)$) to be the maximum number of general
points one can impose on a general $Y \in \caf$.  Since $d$ general points in $\PP^N$ depend
on $dN$ parameters, one has $\ds \mu(\caf)\;\codim(Y,\PP^N) \leq \dim (\caf)$.
For complete intersections one has
\bp \label{ci}
Let $\caf$ be the family of complete intersections of hypersurfaces $Y$ of degrees
$d_1, \ldots, d_r$ in $\PP^N$. Let $s= \min \{d_i\}$. Then
$$ \mu (\caf)= \binom{s+N}{N} - \# \{i \, | \, d_i=s\}.$$
In particular, $\ds \mu(\caf)\;\mbox{\em \codim}(Y,\PP^N)= \dim (\caf)$ if and only if
$d_i=s$ for every $i$.
\ep
We let $\mu(h)$ denote $\mu(\mbox{PGor}(h))$.
There are two obvious upper bounds for $\mu(h)$:
\bp \label{boundmu}
Let $\mu$ denote the maximum number of general points one can impose
on a general $Z\in\mbox{\em PGor}(h)$. Then
\begin{enumerate}
  \item[(a)] $$\mu \leq \frac{1}{3} \dim \mbox{\em PGor}(h).$$
  \item[(b)] if $s=s(h)$, then
$$\mu \leq \binom{s+3}{3} - h^0 (\ideal{Z}(s))=h(s)+ \binom{s+2}{3} $$

\end{enumerate}
\ep

\begin{proof}
For part $(b)$, note that whenever $W \subset Z$, we have $h_Z(n) \geq h_W (n)$.
Hence, if $W$ is a set of $d$ general points contained in
$Z$, we must have
$$h_Z (n) \geq h_W (n) =
\begin{cases}
\begin{array}{lc}
\ds \binom{n+2}{2}  & \mbox{\ \ if $n<s (W)$} \\

\ds d-\binom{s(W)+2}{3}  & \mbox{\ \ if $n=s(W)$.}
\end{array}
\end{cases}
$$
Thus either $s(W) < s $, which means $W$ is contained in a surface of degree $s-1$
and hence $d< \binom{s+2}{3} $, or $s(W)=s$,
in which case $h(s) \geq d-\binom{s+2}{3}$. Thus in any case
$d \leq h(s)+ \binom{s+2}{3}$ as
claimed.

\end{proof}

Lower bounds for $\mu(h)$ are provided by the numbers $\tilde{\mu}$ which we now introduce.
Given a pair $(\tilde{h},m)$ consisting of a $C2$-admissible $h$-vector $\tilde{h}$
and an integer $m$,
we let $h=h(\tilde{h},m)$ denote the $h$-vector determined by the formula
\begin{equation}\label{hh}
\partial h (n)= \tilde{h}(n)-\tilde{h}(m+2-n).
\end{equation}
By Proposition \ref{3.4}, if  $Z$ is a divisor linearly equivalent
to $m H_C-K_C$ on some  ACM curve $C$  with $h$-vector
$\tilde{h}$, then $h=h_Z$. Assuming this linear system is nonempty
for a general curve $C$ in $\mbox{ACM}(\tilde{h})$,
we denote by $U$  the open dense subset of $\mbox{ACM}(\tilde{h})$ consisting of
reduced and locally complete intersection curves, smooth if $\tilde{h}$ is of decreasing type,
on which the linear system $m H_C-K_C$ has the smallest dimension. We can also require
curves in $U$ to have any other general property we may need -- for example the normal bundle
being stable -- as long this does not
make $U$ empty.
Let
$\mathcal{B}=\mathcal{B}(\tilde{h},m)$ be the subscheme of $\mbox{PGor}(h)$ defined as
$$
\mathcal{B}=
\{ Z \in \mbox{PGor}(h) | \, \exists \, C \in U \mbox{ such that } Z \sim m H_C-K_C \mbox{ on } C\}.
$$
By construction $\mathcal{B}$ is irreducible of dimension
\begin{equation}\label{formBB}
\dim \mathcal{B}=  D(\tilde{h},m)  +\dim \mbox{ACM}(\tilde{h})-\delta
\end{equation}
where $D(\tilde{h},m)$ is the dimension of the linear system $|m H_C-K_C|$ for $C \in U$, and
$\delta$ is
the dimension of the family of curves $C \in U$ containing a fixed
general $Z$ in $\mathcal{B}$.

We let $\tilde{\mu}=\tilde{\mu}(\tilde{h},m)=\mu (\mathcal{B})$
denote the maximum number of general points one can impose on a
general $Z \in \mathcal{B}$ . By definition $\mu (h) \geq
\tilde{\mu}(\tilde{h},m)$. We will later see one can often compute
$\tilde{\mu}$.

It is interesting to know when $\mathcal{B}(\tilde{h},m)$ contains a open set of $\mbox{PGor}(h)$,
in which case $\mu=\tilde{\mu}$.
By Theorem \ref{thm25} this is the case when $\tilde{h}$ is the first half of $h$ and
$m \geq 2b(\tilde{h})+1$. It is also the case when $\mathcal{B}(\tilde{h},m)$ has the same dimension
of $\mbox{PGor}(h)$. The latter is known (see Remark \ref{pdim}), while to compute the former
we can use formula (\ref{formBB}). Note that $\dim \mbox{ACM}(\tilde{h})$ is known
(see Proposition \ref{adim}), and,  when the linear system $|m H_C-K_C|$ is nonspecial for a
general $C$
in $\mbox{ACM}(\tilde{h})$, we have $D(\tilde{h},m)= \deg (h)-g(C)$.
Then to determine $\dim \mathcal{B}$ we still need to know the number $\delta$,
although this seems more difficult to compute in general. In several cases
one can show $\delta$ is zero, and of course  this will be the case if the degree of $Z$ is
large enough. We give a precise statement in Corollary \ref{delta-zero} below.

Following Ellia \cite{ellia1}, we denote by $G_{CM} (d,s)$ the maximum genus of an ACM irreducible curve of degree $d$ not lying on a
surface of degree $s-1$. When $\ds d >s(s-1)$, dividing $d$ by $s$ we write $d=st-r$
with $0 \leq r <s$. Then by \cite{GP}
\begin{equation}\label{genus}
G_{CM} (d,s)= 1 +\frac{d}{2} \left( s+\frac{d}{s}-4 \right)-\frac{r(s-r)(s-1)}{2s},
\end{equation}
which is the genus of a curve linked by two surfaces of degrees $s$ and $t$ to a plane curve of degree $r$.

\bp \label{genus-bound}
Assume the base field has characteristic zero. Fix a $C2$-admissible $h$-vector $\tilde{h}$ of decreasing type.
Let $C_1$ and $C_2$ be two distinct irreducible curves in $\mbox{\em ACM} (\tilde{h})$.
Let $d=\deg (C_i)$ and $s=s(C_i)$. If $d \geq 3$, then
$$
g(C_1 \cup C_2) \leq G_{CM} (2d,s).
$$
Furthermore, if equality holds,
then $C_1 \cup C_2$ is an ACM curve linked to plane curve of degree $r$ by two surfaces
of degree $s$ and $t$ respectively, where $r$ and $t$ are defined by the relation $2d=st-r$ with $0 \leq r <s$.
\ep
\begin{proof}
We follow closely the arguments of \cite{ellia1} which require characteristic zero. Note that the numerical character $\chi(X)$ used
by Ellia is an invariant equivalent to the $h$-vector $h_X$ when $X$ is ACM.
Let $C=C_1 \cup C_2$. Since $C_1$ and $C_2$ are the irreducible
components of $C$ and have the same numerical character, the proof of \cite[Theorem 10]{ellia1} shows the
character of $C$ is connected.

Furthermore, since $C_1$ is an ACM curve with $s(C_1)=s$, we have $h_{C_1} (n)=n+1$ for $0 \leq n \leq s-1$, thus
$$
d= \sum h_{C_1}(n) \geq \frac{1}{2}s(s+1).
$$
As in the proof of \cite[Theorem 13]{ellia1}, if we had $g(C)>G_{CM} (2d,s)$, we would have
$\sigma \leq s-1$, where
$\sigma$ is the length of the character $\chi(C)$. Then  $2d \geq s(s+1)> \sigma^2+1$. Since $d \geq 3$,
$C$ contains no curve of degree $2$, hence by \cite[Lemma 12]{ellia1} $C$ is contained in a surface of degree
$\sigma$, which is absurd since $\sigma \leq s-1$.

We conclude $g(C) \leq G_{CM} (2d,s)$. Suppose equality holds and write $2d=st-r$ with $0 \leq r <s$.
Then the argument in  \cite[Theorem 10]{ellia1} together with \cite[Theorem 2.7]{GP} shows
that $C$ is ACM
with the same Hilbert function as a curve linked by two surfaces of degrees $s$ and $t$ to a plane curve of degree $r$.
In particular, $C$ is contained in a surface $S$ of degree $s=s(C_1)$. Therefore $S$ is of minimal degree among surfaces
containing $C_1$, and, since $C_1$ is irreducible, $S$ must also be irreducible. But then looking at the Hilbert function
of $C$ we see $C$ is contained in a surface $T$ of degree $t$ that cuts $S$ properly, and the curve
linked to $C$ by the complete intersection $S \cap T$ is a plane curve of degree $r$.
\end{proof}

\bc \label{delta-zero}
Assume the base field has characteristic zero. Fix an $h$-vector of decreasing type
$\tilde{h}$ and let $U_0$ be the dense open subset of $\mbox{\em ACM} (\tilde{h})$ consisting of smooth irreducible curves.
Let $Z$ be a divisor in the linear system $|m H_C -K_C|$ on some curve $C \in U_0$, and
let $d=\deg(C)$ and $s=s(C)$.

If $d \geq 3$ and
$m d \geq G_{CM} (2 d, s)$,
then $C$ is the unique curve in $U_0$ containing $Z$.
In particular, $\delta$ is zero in this case.
\ec
\begin{proof}
If there was another $C' \in U_0$ containing $Z$, then we would have
$$
\deg (Z)  \leq \deg (C \cap C')= g(C \cup C')-2g(C)+1 \leq G_{CM} (2 d, s)-2g(C)+1.
$$
Since $\deg (Z)= m d-2g(C)+2$, this would imply
$$
m \deg (C) \leq G_{CM} (2 d, s)-1
$$
contradicting the assumptions.
\end{proof}

We now give an example in which $\delta$ is zero, but there is more than one curve
in $\mbox{ACM} (\tilde{h})$ containing $Z$.
Let $C$ be a curve of type $(a,a-1)$ on a smooth quadric surface $Q$. One knows $C$
is ACM with $h$-vector $\tilde{h}=\{1,2, \ldots, 2\}$ ending in degree $b(C)=a-1$
\bl \label{Q}
Suppose $Z$ is a zero dimensional scheme  on a smooth quadric surface that
is the intersection of two  curves $C_1$ and $C_2$ of type $(a,a-1)$ and $(a-1,a)$ respectively. Then $Z$ is AG,
$\deg (Z)= 2a^2-2a+1$ and $m(Z)=2a-3$. If $a \geq 3$ and the curves $C_1$ and $C_2$ are smooth, then these are
the only irreducible curves in $\mbox{\em ACM}(\tilde{h})$ that contain $Z$.
In particular, the dimension $\delta$ of the family of irreducible curves $C \in \mbox{\em ACM}(\tilde{h})$ containing
$Z$ is zero.
\el
\begin{proof}
The fact $Z$ is AG is well known and can be seen as follows: since $Z= C_1 \cap C_2$, we have
$\coo_{C_1} (Z) = \coo_Q (a-1,a) \otimes \coo_{C_1}$, and by the adjunction formula
$Z \sim (2a-3)H-K$ on $C_1$. Therefore $Z$ is AG with $m=2a-3$.

Suppose now $C$ is an irreducible curve in $\mbox{ACM}(\tilde{h})$ which contains $Z$. Then $Z$ is contained in $C \cap C_1$.
By Proposition \ref{genus-bound} and Corollary \ref{delta-zero}, we must have $Z=C \cap C_1$, and $C \cup C_1$ must be the complete
intersection of a quadric surface $Q'$ and a surface of degree $2a-1$.
Since $a \geq 3$, we must have $Q'=Q$, hence
$C$ is a curve type $(a-1,a)$ on $Q$. To finish, observe that there is a one to one correspondence
between curves $C$ of type $(a-1,a)$ on $Q$ and effective divisors linearly equivalent to $Z$ on $C_1$,
because $h^0 \coo_Q(-1,1)=h^1 \coo_Q (-1,1)=0$.
\end{proof}

The following results deal with some particular questions which arise in degree $14$, $21$ and $30$.
\bp \label{14}
Assume the base field has characteristic zero.
The general zero-dimensional arithmetically Gorenstein scheme of degree 14 with $h$-vector
$h=\{ 1,3,6,3,1 \}$ is a divisor
$Z\sim 3H_C-K_C$ on some smooth ACM curve $C$ of degree 6 and genus 3.
\ep
\begin{proof}
The $h$-vector of an ACM curve $C$ of degree $6$ and genus $3$ is $\tilde{h}=\{ 1,2,3 \}$,
the first half of $h$, and $h=h(\tilde{h},3)$, hence a divisor in the linear system
$3H_C-K_C$ has  $h$-vector $\{ 1,3,6,3,1 \}$. Thus it is enough to show
$\dim \mathcal{B}(\tilde{h},3)= \dim \mbox{PGor} (h)$, which is $35$.
By formula (\ref{formBB}) we have
$\dim \mathcal{B}= \dim_C |3H_C-K_C|+\dim \mbox{ACM} (\tilde{h})-\delta=11+24-\delta$,
so it suffices to show $\delta=0$.

Assume by way of contradiction that $\delta >0$: this means that, having fixed a general
$Z$ in $\mathcal{B}$, one can find a positive dimensional family $\mathcal{C}$ of curves $C$ in
$U$ containing $Z$. By the deformation theory of the pair $(Z,C)$ in $\PP^3$, the
infinitesimal deformations of $C$ that leave $Z$ fixed are sections of the normal bundle
$N_C$ that vanish on $Z$. Thus a tangent vector to the  family $\mathcal{C}$ gives a nonzero
section of $N_C(-Z)$. Therefore we obtain a contradiction, and the proposition
will be proven,
if we can show $\mbox{H}^0(C,N_C(-Z))$ is zero.

For this, we claim that, if $C$ is a general ACM curve of degree $6$ and genus $3$
and $Z$ is any divisor of degree $14$ on $C$, then $\mbox{H}^0(C,N_C(-Z))=0$.
Since there are curves of type $(6,3)$ whose normal bundle $N_C$ is stable by \cite{ellia},
and stability is an open property, then for a general such $C$
the normal bundle is stable.
The rank two bundle $N_C$ has degree $28$ while $\deg Z=14$, so we have
$\deg N(-Z)=0$ and, by stability, $\mbox{H}^0(C,N_{C}(-Z))=0$.
\end{proof}

\bp
The general zero-dimensional arithmetically Gorenstein scheme of degree $21$ with $h$-vector
$h=\{1,3,4,5,4,3,1 \}$  is
a divisor
$Z\sim 5 H_C-K_C$ on an ACM curve $C$ of degree 5 and genus 3.
\ep
\begin{proof}
An ACM curve of degree 5 and genus 3 has $h$-vector $\tilde h=\{1,2,1,1\}$.
Since $\dim \mbox{ACM}(\tilde h)=20$ by Corollary \ref{cordimACM},
and $\dim\mbox{PGor}(h)=37$, the statement
follows if $\dim\mbox{PGor}(h)$ equals $\dim\mathcal{B}(\tilde h, 5)=38-\delta$, i.e. if
 $\delta=1$. The general curve of degree 5 and genus 3 is the union of a plane curve $D$ of degree 4 and a line $L$ intersecting $D$ in one point $P$.
 We claim that a general divisor $Z\sim 5 H_C-K_C$  consists of 15 points on the plane quartic and 6 points on the line. In fact, from
 the exact sequences (see \cite[11-10]{sernesi})
\begin{equation}\label{seqser}
\exact{\omega_D}{\omega_C}{\omega_L(P)}
\end{equation}
and
\begin{equation}\label{seqser1}
\exact{\omega_L}{\omega_C}{\omega_D(P)}
\end{equation}
we see $\mathcal{O}_L(5H_C-K_C)=\mathcal{O}_L(6)$ and
$\mathcal{O}_D(5H_C-K_C)=\mathcal{O}_D(4H_D-P)$.
Let $W$ denote the $15$ points of $Z$ lying on the quartic $D$.
Then $W$ is the divisor $4H_D-P$ on $D$, so $h^0(\mathcal{I}_{W,D}(4))=1$ and
$h^0(\mathcal{I}_{W,\mathbb{P}^3}(4))=2$. Therefore $\delta=1$.
\end{proof}

\bp \label{30}
The general zero-dimensional arithmetically Gorenstein scheme $Z$ of degree $30$ with $h$-vector
$h=\{1,3,6,10,6,3,1 \}$  is not of the
form $mH_C-K_C$ on any integral ACM curve.
\ep
\begin{proof}
{From} the $h$-vector we see $m(Z)=5$. We claim that, if $Z\sim 5
H_C-K_C$ on an integral ACM curve $C$, then $C$ has $h$ vector equal to the
first half $k_Z=\{1,2,3,4\}$ of $h_Z$. In fact, by Proposition
\ref{3.4} the $h$-vector of $C$ must be of the form
$h_C=\{1,2,3,4, h_C(4),...., h_C(b)\}$. The curve $C$ then has
degree $d_C = 10+\sum_{n=4}^b h_C(n)=10+\tilde d$ and genus
$g_C=11+\sum_{n=4}^b (n-1) h_C(n)\geq 11 + 3 \tilde d$. Thus
$$
30= \deg Z =
5 d_C -2g_C+2 \leq 50 +5 \tilde d -22 -6 \tilde d +2= 30 -\tilde d,
$$
hence $\tilde d=0$ and $C$ has $h$-vector $\{1,2,3,4\}$.
Now the family of such curves $C$ has dimension $40$, while for fixed $C$,
we have    $\dim_C |Z|=h^0\mathcal{O}_C(5H_C-K_C)-1=19$  because
$\deg Z > 2g_C-2$ (here we use the hypothesis $C$ integral).
Thus the family of schemes $Z$ of the form
$5H_C-K_C$ on some integral ACM curve $C$ has
dimension at most 59. However, if we compute  $\dim\mbox{PGor}(h_Z)$,
as explained in Remark \ref{pdim} we find $\dim\mbox{PGor}(h_Z)=63$. Thus a general
$Z$ cannot be of the form $mH_C-K_C$ on any integral ACM curve $C$.
\end{proof}

\br It seems unlikely that $Z$ could be of the form $-K+mH$ on any ACM curve (possibly
reducible), but we do not have a complete proof.
\er

In Table \ref{tavola3} we list, for every degree $d \leq 30$,
all possible $h$-vectors of zero-schemes  of degree $d$ in  $\mathbb{P}^3$ not contained in a plane.
The list is constructed using Proposition \ref{G3}.
For every $h$ vector in the table we record
\begin{itemize}
  \item the dimension $A$ of $\mbox{PGor}(h)$, which can be computed as explained in Remark \ref{pdim},
  or applying the formula of \cite{K};
  \item the invariant $m$ of Propositions \ref{bb}  and \ref{bb1};
  \item at least one $h$-vector $\tilde{h}$ such that $h=h(\tilde{h},m)$ as in formula (\ref{hh});
  \item the degree and genus of an ACM curve with $h$-vector $\tilde{h}$;
  \item the dimension $B$ of $\mathcal{B}(\tilde{h},m)$; we note that:
\begin{enumerate}
  \item[(a)] when $h$ is the $h$-vector of a complete intersection of type $(a,b,c)$ and $\tilde{h}$ is
   the $h$-vector of a complete intersection of type $(a,b)$, then $B=A$.
  \item[(b)] when Theorem \ref{thm25}(c) holds, that is, $m \geq 2 b(\tilde{h})+1$, one has $B=A$.
  \item[(c)] in all other cases, since in the table we always have $\deg (h) \geq 2g(\tilde{h})-1$,
  formula~(\ref{formBB}) gives
$B=  \deg (h)-g(\tilde{h})  +\dim \mbox{ACM}(\tilde{h})-\delta$;
\end{enumerate}
  \item the maximum number $\tilde{\nu}$ of general points on a general curve $C \in \mbox{ACM}(\tilde{h})$; by results
  of Perrin \cite{P} and Ellia \cite{ellia1}, it happens that in all cases in our table one has
  $$
 \tilde{\nu} = \mbox{Min} \left( \frac{1}{2} \dim \mbox{ACM}(\tilde{h}), \; \alpha \, \right)
  $$
where $\alpha$ is the dimension of the family of surfaces of degree $s(C)$ that contain curves in $\mbox{ACM}(\tilde{h})$.
When $s(C) \leq 3$, we have $\alpha={{s(C)+3}\choose 3}-1$.
  Note that the problem of determining $\tilde \nu$ for ACM curves of higher degree remains open.
  \item the maximum number $\tilde{\mu}$ of general points on a general $Z  \in \mathcal{B}(\tilde{h},m)$;
  for $C$ a general curve in $\mbox{ACM}(\tilde{h})$, we have $\tilde{\mu}= \tilde{\nu}$ if
  $\dim |m H_C-K_C| \geq  \tilde{\nu}$, while, if $\dim |m H_C-K_C| <  \tilde{\nu}$, we can only say
  $$\dim |m H_C-K_C| \leq \tilde{\mu} \leq \mu;$$
    \item the maximum number $\mu$ of general points on a general $Z  \in \mbox{PGor}(h)$;  note that $\mu=\tilde{\mu}$
  when $A=B$.
\end{itemize}
We indicate with a check ''$\checkmark$" whether a general AG zero-scheme is in the linear system $|mH-K|$ on some
ACM curve $C$ with $h$-vector $\tilde{h}$.


\section{General sets of points in $\mathbb{P}^3$}
\bp \label{glicci} A set of $n\leq 19$ general points in $\mathbb{P}^3$ can be obtained by a sequence of ascending (i.e.
degree increasing) Gorenstein liaisons from a point. In
particular, it is glicci. \ep
\begin{proof} If we have a set $W$ of $n$ general points and we can find a family of AG schemes
$Z$ of degree $d$ containing $\mu$ general points, and if
$\frac{1}{2}d<n\leq \mu$, then we can perform a descending
Gorenstein liaison from $W$ using $Z$ to get a new set $W'$ of
$d-n<\frac{1}{2} d<n$ points. Since the process is reversible, the
set $W'$ also consists of general points, and we can
continue the process. For one or two points the result is trivial.
For $n\geq 3$ general points, using Table 1, we choose
$Z$ of degree $d$ and $h$-vector $h$ as follows
$$
\begin{array}{llll}
n & d \ \ \ \ & h & \mu \\
3,4 & 5& \{ 1,3,1\} & 5 \\
5,6,7 & 8& \{ 1,3,3,1\} & 7 \\
8,9,10,11 \ \ & 14& \{ 1,3,6,3,1\} & 11 \\
12,13,14 & 20 & \{ 1,3,6,6,3,1\} & 14 \\
15,16,17 & 27 & \{ 1,3,6,7,6,3,1\} & 17 \\
18,19 & 30 & \{ 1,3,6,10,6,3,1\} \ \ & \geq 19 \\
\end{array}
$$
\end{proof}
At present, this is as far as we can go, because for $n=20$ we do
not know if an AG zero scheme of degree $30$ and $h=\{
1,3,6,10,6,3,1\}$ has $\mu\geq 20$.

\bt \label{dl}
A set of $n\geq 56$ general points in $\PP^3$ admits no strictly
descending Gorenstein liaison.
\et
\begin{proof}
Let $W$ be a set of $n$ general points, with $s(W)=s$, so that its $h$-vector is
$$
h_W=\left\{ 1,3,6,\ldots , \binom{s+1}{2},a\right\}
$$
with $0\leq a<\binom{s+2}{2}$ and $n=\binom{s+2}{3}+a$. If $W$ is contained in an
AG scheme $Z$ then $s(Z)\geq s$. On the other hand, if the residual scheme $W'$ has degree
less than $n$, then looking at the $h$-vectors and using Proposition \ref{variation},
we see that there are only three possibilities for
$h_Z$:

\vspace{.2truecm}
\begin{enumerate}
\item[Type 1] \ \ $h_Z=\{1,3,6,\ldots, {{s+1}\choose 2}, \ldots, 6,3,1\}$
\item[Type 2] \ \ $h_Z=\{1,3,6,\ldots, {{s+1}\choose 2},{{s+1}\choose 2}, \ldots, 6,3,1\}$ and
$a > 0$
\item[Type 3] \ \ $h_Z=\{1,3,6,\ldots, {{s+1}\choose 2},b,{{s+1}\choose 2}, \ldots, 6,3,1\}$,
with $\binom{s+1}{2}\leq b\leq\binom{s+2}{2}$ and $a > \displaystyle\frac{1}{2}b$
\vspace{.2truecm}
\end{enumerate}
A necessary condition for $Z$ to contain $n$ general points is that $\dim\mbox{PGor}(h)\geq 3n$.
We compute $\dim\mbox{PGor}(h)$ for each of the types above, using Corollary \ref{dimension}
and induction on $s$. Setting $b={{s+1}\choose 2}+c$, so that $0\leq c\leq s+1$, we find
\vspace{.2truecm}
\begin{enumerate}
\item[Type 1]  \ \ $\dim\mbox{PGor}(h)=4s^2-1$
\item[Type 2]  \ \ $\dim\mbox{PGor}(h)=4s^2+3s-1$
\item[Type 3]  \ \ $\dim\mbox{PGor}(h)=4s^2+4s+4c-1$.
\end{enumerate}
\vspace{.2truecm}
Now writing the inequality $\dim\mbox{PGor}(h)\geq 3n$, we find for Type 1
$$
4s^2-1 \geq 3\binom{s+2}{3}+3a,
$$
and using $a\geq 0$, this implies $s<5$.
\par\noindent
For Type 2, we find
$$
4s^2+3s-1\geq 3\binom{s+2}{3}+3a.
$$
Again using $a>0$, this implies $s<6$.
\par\noindent
For Type 3, we have
$$
4s^2+4s+4c-1\geq 3\binom{s+2}{3}+3a.
$$
Using $a > \displaystyle\frac{1}{2} b=\displaystyle\frac{1}{2} \binom{s+1}{2}+\displaystyle\frac{1}{2}c$ this gives
$$
4s^2+4s+\frac{5}{2} c-1\geq 3\binom{s+2}{3}+\frac{3}{2}\binom{s+1}{2}.
$$
Now using $c\leq s+1$ we get
$$
4s^2+4s+\frac{5}{2} (s+1) -1\geq 3\binom{s+2}{3}+\frac{3}{2}\binom{s+1}{2}
$$
which implies $s<6$

Thus for $s\geq 6$, and hence for  $n\geq 56$, a set of $n$ general points has no
descending Gorenstein liaison.
\end{proof}

\br
Checking possible values of $a$ and $c$ for $s=5$, the same method applies
for all $n\geq 35$, except possibly $36$, $37$, $38$, $45$, $46$, $47$.
\er

One of us \cite[Proposition 2.7]{H3} has shown that a set of $n < 19$, $n \neq 17$, general points in $\PP^3$ can be obtained
by a sequence of ascending biliaisons from a point. On the other hand, we can prove:
\bt \label{db}
A set of $n\geq 56$ general points in $\PP^3$ admits no strictly
descending elementary biliaison.
\et
\begin{proof}
Suppose a set $Z$ of  $n\geq 56$ general points  admits a descending biliaison on an ACM curve $C$, i.e., $Z \sim W+H$ on $C$.
Then by Proposition \ref{variation}
the $h$-vectors satisfy
$$ h_Z(l)=h_C (l) +h_W(l+1)$$
for all $l$. Let $s=s(Z)$, so that $h_Z (s-1)=\frac{1}{2}s(s+1)$. It follows that $h_C (s-1)$ and $h_W(s-2)$ must each achieve
their maximum values, namely $s$ and $\frac{1}{2}s(s-1)$ respectively.
It follows also that $h_C(l)=0$ for $l \geq s+1$, since this is
true also for $h_Z$ by formula (\ref{hgeneral}).
Thus $C$ satisfies the hypothesis of Corollary \ref{cordimACM}
and the dimension of the family $\mbox{ACM} (h_C)$ is $4d$, where $d$ is the degree of $C$.

Now in order for $C$ to contain $n$ general points, we must have $\dim \mbox{ACM} (h_C) \geq 2n$.
Let $h_C(s)=a$. Then $d= \deg \, C= \frac{1}{2}s(s+1)+a$. Let $h_Z (s)= b$. Then
$n = \deg \, Z= \binom{s+2}{3}+b$. Furthermore, $a \leq s+1$ and $a \leq b$.

{From} $4d \geq 2n$, we thus obtain $2d \geq n$, or
$$
s(s+1)+2a \geq \binom{s+2}{3}+b
$$
Writing $2a \leq s+1+b$, we get $$(s+1)^2 \geq \binom{s+2}{3}$$
which implies $s \leq 2 + \sqrt{10} < 6$. So for $s \geq 6$,
and hence for any $n \geq 56$, a set $Z$ of $n$ general points admits no descending biliaison.
\end{proof}

\br
The same argument, taking into account the exact values of $a$ and $b$, applies to all $n \geq 31$,
except for $n= 40, 41, 42$. Using a slightly more sophisticated argument,
we can treat smaller values of $n$ such as the following case of $n=20$.
\er

\bex
A set $Z$ of $20$ general points admits no descending elementary biliaison. Indeed,
by the analysis in the proof above, the only possibility would be on an ACM curve $C$ with $h$-vector
$h_C=\{ 1,2,3,4 \}$. This is a curve of degree $10$ and genus $11$. Since the family of all $Z$'s has dimension $60$, and
the family of pairs $(C,Z)$ with $Z \subset C$ has dimension $60$ also,
we conclude that a general $Z$ must lie on a general $C$, and that the points in $Z$ are also general on $C$.
Thus for a general $Z$, the divisor $W=Z-H$ on $C$ is a general divisor of degree $10$. But since the genus of $C$
is $11$, this general $W$ cannot be effective. Thus the general $Z$ has no descending biliaison.

On the other hand, we do not know if $Z$ admits a strictly descending Gorenstein liaison,
because we cannot answer the question whether there are AG zero-schemes of degree $30$ containing $20$ general points.
If so, we could link $20$ general
points to $10$ general points and thus $20$ general points would be
glicci. In our notation, the question is whether $\mu\geq 20$. We have
only been able to show $19\leq \mu \leq 21$.
There are AG schemes of degree
$30$ in the linear system $5H-K$ on an ACM curve with $(d,g)=(10,11)$, but
these can contain at most $19$ general points. The more general
AG schemes of degree $30$ are not of the form $mH-K$ on any integral ACM curve by Proposition \ref{30}, so
some new technique will be necessary to answer this question.
\eex

\bex We can show by an analogous but more complicated argument
that a set $Z$ of $31$ general points on a nonsingular cubic surface in
$\mathbb{P}^3$ does not admit any descending Gorenstein biliaison
in $\mathbb{P}^3$. Since $Z$ is in the strict Gorenstein liaison
equivalence class of a point \cite[2.4]{H3}, it is glicci.
Furthermore, since in codimension 3, even strict Gorenstein
liaison gives the same equivalence relation as Gorenstein
biliaison \cite[5.1]{H1}, $Z$ is even in the Gorenstein biliaison
equivalence class of a point. This is the first example we know of
of a scheme $Z$ that is glicci but cannot be obtained by a
sequence of ascending Gorenstein biliaisons from a linear space.
\eex

\section{Conclusion}
We have established a number of fundamental results about arithmetically Gorenstein zero-dimensional
schemes in $\mathbb{P}^3$. In particular, we have investigated those that occur in the form
$mH-K$ on an ACM curve and we studied the number of general points that one can impose on an AG scheme with given
$h$-vector, in order to understand the possible Gorenstein liaisons that one can perform on
a set of general points.
In all cases we are aware of where a class of zero-dimensional subschemes of
$\mathbb{P}^3$ has been proved to be glicci, the proof was actually accomplished using strict
Gorenstein liaisons, i.e. using only  those AG schemes of the form $mH-K$ on some ACM curve (see \cite[$\S$ 1]{H3},
for the terminology of strict G-liaison).
Remembering that a Gorenstein biliaison is a composition of
two strict G-liaisons, this remark applies to the determinantal schemes of \cite[3.6]{KMMNP},
to any zero-scheme on a non singular quadric surface or a quadric cone \cite[5.1 and 6.1]{CH}
and to $n$ general points on a non
singular cubic surface in $\mathbb{P}^3$  \cite[2.4]{H3}. For us, this underlines
the importance of studying those
more general AG zero-dimensional schemes not of the form $mH-K$ on any ACM curve, and by making use
of them either to prove or disprove the assertion that ''Every zero-scheme in $\mathbb{P}^3$ is
glicci''.

In the course of this work we have been led to reconsider some old problems whose solution would
be helpful in furthering the work of this paper. One is the problem of Perrin's thesis \cite{P}
to find how many general points one can impose to a curve of given degree and genus in $\mathbb{P}^3$.
For us it is the ACM curves that are relevant, so we ask: is it true that for a general smooth
ACM curve with $h$-vector $\tilde h$, the number $\tilde\nu$ of general points
one can impose on the curve is given by the formula mentioned in section 6:
$$
\tilde\nu=\min \left( \frac{1}{2}\dim\mbox{ACM}(\tilde h), \alpha \right) ?
$$
The other old question, which appeared in a special case in the proof of Proposition \ref{14},
concerns the stability of the normal bundle of a space curve. The problem was stated
in \cite{HC}, and has been more recently studied by Ellia \cite{ellia1}: if $C$ is a general smooth ACM
curve of degree $d$, genus $g$, $s=s(C)$ and
$$
g< d(s-2)+1 \qquad ({\rm resp.}\ \leq)
$$
then is the normal bundle of $C$ stable (resp. semistable)?

\vfill
\noindent
\hoffset -0.6truecm
\begin{table}[p]
\caption{Nondegenerate AG zero-dimensional subschemes in $\mathbb{P}^3$ of degree $\leq 30$ (char. $k = 0$)}
\label{tavola3}
\footnotesize{
\begin{tabular}{|c|l|c|l|c|c|c|c|c|c|l|}
\hline \hline
&&&&&&&&&&\\
$d$ & \;\;\;\;\;\;\;$h$  &  $A$ & $m$ & $\tilde{h}$ &$(\tilde{d},\tilde{g})$&  $B$  & $\ $&$\tilde{\nu}$
& $\tilde{\mu}$   & $\mu$  \\
&&&&&&&&&& \\
\hline
5&$\{ 1,3,1 \}$           & 15     & 1 &$\{1,2\}$      & (3,0)   & 15      &$\checkmark$  &  6 & 5  &$5^a$
\\
\hline
8 & $\{ 1,3,3,1 \}$       & 21     & 2 &$\{1,2\}$      & (3,0)   &  20     & no           &  6 & 6 & 7$^e$
\\
\ &                   \   & \      &   &$\{1,2,1\}$    & (4,1)   &  21$^e$ & $\checkmark$ &  8 & 7 &
\\
\hline
11 & $\{ 1,3,3,3,1 \}$    & 23     & 3 & $\{ 1,2\}$    & (3,0)   &  23$^d$ & $\checkmark$ &  6 & 6 & $6^d$
\\
\hline
12 &  $\{ 1,3,4,3,1 \}$   & 27     & 3 & $\{ 1,2,1\}$  & (4,1)   & 27$^e$ & $\checkmark$  &  8 & 8 & $8^e$
\\
\hline
13 & $\{ 1,3,5,3,1 \}$    &  31    & 3 & $\{ 1,2,2\}$  & (5,2)   &      31&  $\checkmark$ &  9 & 9 & $9^b$
\\
\hline
14 & $\{ 1,3,3,3,3,1 \}$  &  26    &4  & $\{ 1,2\}$    & (3,0)   &  26$^d$ & $\checkmark$ &  6 & 6 & $6^d$
\\
\hline
14 & $\{ 1,3,6,3,1 \}$    &  35    & 3 & $\{ 1,2,3\}$  & (6,3)   &    35   & $\checkmark$   & 12& 11& $11^a$
\\
\hline
16 & $\{ 1,3,4,4,3,1 \}$  & 31     & 4 & $\{ 1,2,1\}$  & (4,1)   &  31$^e$ & $\checkmark$   & 8 & 8 & $8^e$
\\
\hline
17 & $\{1,3,3,3,3,3,1 \}$ & 29     & 5 & $\{ 1,2\}$    & (3,0)   &  29$^d$ & $\checkmark$   & 6 & 6 & $6^d$
\\
\hline
18 & $\{ 1,3,5,5,3,1 \}$  &  37    & 4 & $\{ 1,2,2\}$  & (5,2)   &      36 & no             & 9 & 9 & $9^e$
\\
\  &                  \   &        & \ & $\{1,2,2,1\}$ & (6,4)   &  37$^e$ & $\checkmark$   & 9 & 9 & \
\\
\hline
20 &$\{1,3,3,3,3,3,3,1 \}$&  32    & 6 & $\{ 1,2\}$    & (3,0)   & 32$^d$  & $\checkmark$   & 6 & 6& $6^d$
\\
\hline
20 & $\{ 1,3,4,4,4,3,1 \}$&  35    & 5 & $\{ 1,2,1\}$  & (4,1)   & 35$^e$  & $\checkmark$  & 8 & 8 &$8^e$
\\
\hline
20 &$\{ 1,3,6,6,3,1 \}$   &  44    & 4 & $\{ 1,2,3\}$  & (6,3)   &      41 & no            & 12&12 &$14^a$
\\
\  & \                    &  \     &\  &$\{1,2,3,1\}$ & (7,5)    &      43 & no            & 14&14 & \
\\
\  &\                     &   \    &\  &$\{1,2,3,2\}$ & (8,7)    & $\leq 44$& ?            & 16&  ?&\
\\
\hline
21 & $\{ 1,3,4,5,4,3,1 \}$& 37     & 5 &$\{ 1,2,1,1\}$& (5,3)    &    37     & $\checkmark$ &  4 & 4  &$4$
\\
\hline
23&$\{1,3,3,3,3,3,3,3,1\}$&  35    & 7 & $\{ 1,2\}$  & (3,0)     & 35$^d$   & $\checkmark$&  6 & 6  & $6^d$
\\
\hline
23& $\{1,3,5,5,5,3,1 \}$  &  41    & 5 & $\{ 1,2,2\}$& (5,2)     & 41$^d$   & $\checkmark$  & 9 & 9  & 9$^d$
\\
\hline
24& $\{ 1,3,4,4,4,4,3,1\}$&  39    & 6 & $\{ 1,2,1\}$& (4,1)     & 39$^d$   & $\checkmark$  & 8 & 8  & $8^{d}$
\\
\hline
24 & $\{ 1,3,5,6,5,3,1 \}$&  44    & 5 &$\{1,2,2,1\}$& (6,4)     & 44$^e$   & $\checkmark$  & 9 & 9  & $9^{e}$
\\
\hline
25 &$\{1,3,5,7,5,3,1 \}$  &  47    & 5 & $\{1,2,2,2\}$& (7,6)    & 47       & $\checkmark$  & 9 & 9 & $9^{a,b}$
\\
\hline
26&$\{1,3,3,3,3,3,3,3,3,1\}$& 38   & 8 & $\{ 1,2\}$  & (3,0)     & 38$^d$   & $\checkmark$  & 6 & 6   & $6^d$
\\
\hline
26&$\{1,3,4,5,5,4,3,1\}$  &  43    & 6 & $\{1,2,1,1\}$& (5,3)    & 43       & $\checkmark$  &4  & 4 & $4$
\\
\hline
26&$\{1,3,6,6,6,3,1\}$    &  47    & 5 & $\{ 1,2,3\}$& (6,3)     & 47$^d$   & $\checkmark$  &12 & 12  & $12^d$
\\
\hline
27&$\{ 1,3,6,7,6,3,1 \}$  &  51    & 5 & $\{1,2,3,1\}$& (7,5)    & 50       & no            &14 & 14& $17^e$
\\
\ & \                     & \      &\  &$\{1,2,3,2,1\}$& (9,10)  & 51$^e$   & $\checkmark$  &18 & 17& \
\\
\hline
28&$\{1,3,4,4,4,4,4,3,1\}$& 43     & 7 & $\{ 1,2,1\}$ & (4,1)    & 43$^d$   & $\checkmark$  & 8 &  8& $8^{d}$
\\
\hline
28 &$\{ 1,3,5,5,5,5,3,1\}$& 46     & 6 & $\{ 1,2,2\}$ & (5,2)    & 46$^d$   & $\checkmark$  & 9 & 9 & $9^d$
\\
\hline
28  & $\{ 1,3,6,8,6,3,1 \}$& 55    & 5 &$\{ 1,2,3,2\}$& (8,7)    &       53 & no            & 16& 16 & 16$\leq\mu\leq$18
\\
\   &                    \ &\       &\  &$\{1,2,3,3,1\}$&(10,12) & $\leq55$ & ? & \         19& 16  &
\\
\hline
29&$\{ 1,3,3,3,3,3,3,3,3,3,1 \}$ & 41  & 9 & $\{ 1,2\}$ & (3,0)  & 41$^d$ & $\checkmark$&6 & 6
 & $6^d$
\\
\hline
29  & $\{ 1,3,6,9,6,3,1 \}$& 59  & 5 & $\{ 1,2,3,3\}$ & (9,9) &$\leq 56$  &  no  & 18& 18
& 18$\leq\mu\leq$19 \\
   &                 &   \ &\     &$\{1,2,3,4,1\}$&(11,14) &$\leq 59$  & ?   & $22$& ?  &
\\
\hline
30  & $\{ 1,3,5,6,6,5,3,1 \}$ & 50  & 6 & $\{ 1,2,2,1\}$ & (6,4) &50$^e$ & $\checkmark$ & 9& 9
 & 9$^{e}$
\\
\hline 30  & $\{ 1,3,6,10,6,3,1 \}$ & 63 & 5 & $\{
1,2,3,4\}$ & (10,11) & $\leq 59^f$  & ? &20 & 19    & 19$\leq\mu\leq $21
\\
\hline
\multicolumn{10}{l}{$a$: $\mu=\tilde{\mu}$ because $A=B$ or $\tilde{\mu}=[A/3]$}\\
\multicolumn{10}{l}{$b$: the upper bound on $\mu$ is given by Proposition \ref{boundmu}}\\
\multicolumn{10}{l}{$c$: $\delta=0$ by \ref{delta-zero}, \ref{Q} or \ref{14} } \\
\multicolumn{10}{l}{$d$:  Theorem \ref{thm25}(c) applies, and $\tilde{\mu}= \tilde{\nu}$ because
  $\dim |m H_C-K_C| \geq  \tilde{\nu}$}
\\
\multicolumn{10}{l}{$e$:  complete intersection case}\\
\multicolumn{10}{l}{$f$: see Proposition \ref{30}}
\\
\end{tabular}
}
\end{table}

\end{document}